\documentclass[11pt]{article}
\usepackage{graphicx}
\usepackage{enumerate}
\usepackage{folded}

\title{Pseudo-holomorphic maps into folded symplectic four-manifolds}
\author{Jens von Bergmann}

\begin{document}

\thispagestyle{empty}
{
  \maketitle 
  
  \abstract{Every oriented 4-manifold admits a folded symplectic
    structure, which in turn determines a homotopy class of compatible
    almost complex structures that are discontinuous across the
    folding hypersurface (``fold'') in a controlled fashion. We define
    {\em folded holomorphic maps}, i.e. pseudo-holomorphic maps that
    are discontinuous across the fold. The boundary values on the fold
    are mediated by {\em tunneling maps} which are punctured
    $\mathcal{H}$-holomorphic maps into the folding hypersurface with
    prescribed asymptotics on closed characteristics.
  
    Our main result is that the linearized operator of this boundary
    value problem is Fredholm and thus we obtain well behaved local
    finite dimensional moduli spaces.
  
    As examples we characterize the moduli space of maps into folded
    elliptic fibration $E^F(1)$ and we construct examples of degree
    $d$ rational maps into $S^4$. Moreover we explicitly give the
    moduli space of degree 1 rational maps into $S^4$ and show that it
    possesses a natural compactification.
    
    This aims to generalize the tools of holomorphic maps to all
    oriented 4-manifolds by utilizing folded symplectic structures as
    opposed to other types of pre-symplectic structures as done in
    \cite{taubes_beasts}.
}

  \tableofcontents
}

\section{Introduction}
\label{sec:overview}
In his seminal 1984 paper \cite{gromov}, M. Gromov showed how to
extend tools of complex geometry to the symplectic category. In the
last decade this has led to a vibrant new field based on the study of
``pseudo-holomorphic curves'' in symplectic manifolds, yielding many
powerful results regarding invariants and the recognition of
symplectic manifolds. 

Unfortunately, these methods do not apply to all manifolds -- many
manifolds do not admit symplectic forms. There are different
possibilities to extend the theory of pseudo-holomorphic curves to a
broader class of manifolds.  One approach, being pursued by C.
Taubes, begins with the observation that every compact oriented
4-manifold with intersection form that is not negative definite admits
a closed 2-form that degenerates along a disjoint union of circles.
Taubes has made a detailed study of the behavior of pseudo-holomorphic
curves approaching these circles (\cite{taubes_beasts}).

In \cite{ana1} A. Cannas, V. Guillemin and C. Woodward introduced the
notion of folded symplectic structures, which we describe in Section
\ref{sec:folded_symplectic_manifolds}. Every orientable 4-manifold
admits a folded symplectic structure (\cite{ana2}).

In this article we construct finite-dimensional moduli spaces of
folded holomorphic maps into closed oriented folded symplectic
4-manifolds with circle invariant folds.

The 4-sphere $S^4$ with its canonical folded symplectic structure is
of this form. We compute the moduli space of folded holomorphic maps
of degree 1 into $S^4$ and give examples of higher degree maps.

Folded holomorphic maps are discontinuous across the fold, a feature
that is in stark contrast to properties of pseudo-holomorphic maps. We
explain why this is natural and also necessary in order to obtain a
well defined moduli space. The discontinuity is controlled by
tunneling maps which are the central objects of this
theory. Intuitively, folded holomorphic maps can be thought of
``smoothings'' of pseudo-holomorphic maps into each component of the
complement of the folding hypersurface equipped with degenerate
structures as in Symplectic Field Theory \cite{SFT} or in \cite{IP_rel}.

\vskip1cm
\subsection{Folded Symplectic Manifolds}
\label{sec:folded_symplectic_manifolds}
We begin by recalling some definitions and basic properties of folded
symplectic manifolds.

\begin{definition}[Folded Symplectic Structure]
  Let $X$ be a smooth $2n$-dimensional manifold. A {\em folded
    symplectic structure} $\omega$ is a closed 2-form such that
  $\omega^{n}$ is transverse to 0 (so $Z=(\omega^n)^{-1}(0)$ is a
  smooth codimension 1 hypersurface) and $\omega^{n-1}|_Z$ is
  non-vanishing. $Z$ is called the {\em fold}.
\end{definition}

The last condition means that the kernel $E=\ker(\omega)$ of $\omega$,
which by transversality is a real 2-plane bundle over $Z$, is
transverse to $TZ$. This is equivalent to the requirement that
\begin{eqnarray*}
  L=\ker(\omega|_Z)\subset TZ 
\end{eqnarray*}
is a 1-dimensional foliation. $L$ is called the {\sl characteristic
  foliation}.
Thus intrinsically, the fold of a folded symplectic manifold is
indistinguishable from an orientable hypersurface in a symplectic
manifold.

Folded symplectic structure, like symplectic and contact structures,
can locally be put in standard form.
\begin{theorem}[Darboux]
  For every folded symplectic form $\omega$ there exist local
  coordinates near the fold such that $\omega$ has the form
  \begin{eqnarray*}
    x_1\,dx_1\wedge dx_2+dx_3\wedge dx_4,
  \end{eqnarray*}
  where the fold is locally given by $\{x_1=0\}$.
\end{theorem}

More generally, in \cite{ana1} it is proved that for any
$\alpha\in\Omega^1(Z)$ that does not vanish on $L$ we can extend the
inclusion $i:Z\hookrightarrow X$ of the fold to an orientation
preserving diffeomorphism
\begin{eqnarray}
  \phi:(-\varepsilon,\varepsilon)\times Z\rightarrow U
  \label{eq:darboux_coords}
\end{eqnarray}
onto a tubular neighborhood $U$ of the fold such that
\begin{eqnarray}
  \phi^\ast\omega&=&\pi^\ast
  i^\ast\omega+\frac12 d(r^2\pi^\ast\alpha)\label{eq:darboux},
\end{eqnarray}
where $r$ is the coordinate function on $(-\varepsilon,\varepsilon)$
and $\pi$ is the projection $\pi(r,z)\rightarrow z$. By possibly
changing the volume form we assume for simplicity that $r=det(\omega)$
on $U$.

\begin{definition}
  A {\em morphism} $\psi:X_1\rightarrow X_2$ of folded symplectic
  manifolds $(X_1,\omega_1)$ and $(X_2,\omega_2)$ is a diffeomorphism
  satisfying
  \begin{eqnarray*}
    \psi^\ast\omega_2=\omega_1.
  \end{eqnarray*}
\end{definition}
Such morphisms automatically take folds to folds.

It was shown in \cite{ana2} that every oriented 4-manifold admits a
folded symplectic structure. To define folded holomorphic maps we
require a special structure along the fold. The following Definition
borrows from \cite{hofer_zehnder} and \cite{SFT_compactness}.
\begin{definition}
  A folded symplectic manifold $(Z,\omega)$ together with a choice of
  1 form $\alpha$ on the fold $Z$ is called {\em stable} if
  \begin{eqnarray}
    \alpha\wedge\omega>0\qquad\mathrm{and}\quad
    \ker(\omega|_Z)\subset\ker(d\alpha).\label{eq:stable}
  \end{eqnarray}
\end{definition}
Such a choice of 1-form $\alpha$ defines an almost contact form
$(\omega,\alpha)$ on $Z$. In particular $\alpha$ defines a section $R$
of the line bundle $L$ by the normalizing condition
\begin{eqnarray*}
  \alpha(R)=1.
\end{eqnarray*}
The flow generated by $R$ induces an $\setR$ action on $Z$ preserving
$\omega$ and $\alpha$. The kernel of $\alpha$ gives a symplectic
subbundle $(F=\ker(\alpha),\omega|_F)$ over $Z$ such that
\begin{eqnarray}\label{eq:E-F-splitting}
  T_ZX=E\oplus F.
\end{eqnarray}

A {\em parametrized closed characteristic} is an injective map
$x:S^1\rightarrow Z$ satisfying
\begin{eqnarray*}
  dx(\del_\theta)=\frac{1}{\tau}\cdot R,\ 
  \mathrm{for\ some\ }\tau\in(0,\infty).
\end{eqnarray*}
A {\em closed characteristic} is an image of a parametrized closed
characteristic and we let $\mathcal R$ be the set of closed
characteristics equipped with the disjoint union topology. For each
closed characteristic $x\in\mathcal{R}$, the corresponding constant
$\tau=\tau_x$ for a parametrized closed characteristic representing it
is called the {\em minimal period}. In the case of $S^1$-invariant
folds we have that $\tau_x=1$ for all $x\in\mathcal{R}$.

Now choose a background metric $g'$ on $X$. To simplify future
computations we assume that in the Darboux chart (\ref{eq:darboux})
\begin{eqnarray*}
  \phi^\ast g'|_Z=dr\otimes dr
  +\alpha\otimes\alpha+g_F
\end{eqnarray*}
where $g_F$ is a metric on $F$ that is compatible with $\omega|_Z$.
Mimicking the standard procedure to generate a compatible triple on a
symplectic manifold using the background metric $g'$ we obtain a
folded triple $(\omega,g,J)$ on $X\setminus Z$ satisfying the
compatibility conditions
\begin{eqnarray*}
  J^\ast\omega&=&\omega\\
  g(u,v)&=&\omega(u,Jv).
\end{eqnarray*}
Near the fold the triple has the following ``standard'' form.
\begin{lemma}\label{lem:triple}
  The complex structure $J|_{X_\pm\setminus Z}$ has one-sided limits
  $J_\pm$ on $T_ZX$ such that the splitting (\ref{eq:E-F-splitting})
  is $J_\pm$-invariant and $J_+|_F=J_-|_F$ and $J_+|_E=-J_-|_E$.
  
  Moreover, $E$ and $F$ extend to a transverse pair of $J$-invariant
  subbundles of $TX$ to a neighborhood $U\subset X$ of $Z$, also
  denoted by $E$ and $F$, such that
  \begin{eqnarray}\label{eq:triple_on_E}
    \omega_E=\det(\omega)\mu,\qquad
    g_E=|\det(\omega)|h,\qquad 
    J_E=\sign{\det(\omega)}\tilde J
  \end{eqnarray}
  where $(\mu,h,\tilde J)$ is a smooth compatible triple on $E$ with
  $\tilde J\del_r=R$ on $Z$.
\end{lemma}
\begin{proof}
  Define the $g'$-skew-endomorphism $A$ of $TX$ by $\omega(u,v)=g'(Au,v).$
  This exists and is unique since $g'$ is positive definite. Observe
  that over $Z$, $A|_{E}=0$ and $A|_{F}$ is non-degenerate.  Since $Z$
  is compact there exists $\delta>0$ with $|\lambda|>\delta$ for all
  eigenvalues of $A|_{F}$.  Thus there exists a neighborhood $U$
  of $Z$ such that $|\lambda_1|<\delta/2<|\lambda|$, over $U$, where
  $\pm\lambda_1$ are the smallest eigenvalues of $A$ and $\pm\lambda$
  are the second smallest eigenvalues of $A$.  Define $E$ to be the
  2-plane bundle over $U$ given by the real eigenspace spanned by
  $\pm\lambda_1$ and set $F=E^{\perp_{g'}}$. Then for $u\in E$ and
  $v\in F$ we have $Au=\tilde u\in E$ and
  \begin{eqnarray*}
    \omega(u,v)=g'(Au,v)=g'(\tilde u,v)=0.    
  \end{eqnarray*}
  Therefore $\omega$ splits as $\omega_E\oplus \omega_{F}$. $A$ also
  splits as $A=A_E\oplus A_{F}$ because the real eigenspaces
  corresponding to different pairs of complex conjugate eigenvalues of
  a skew-endomorphism are orthogonal. Since $A_{F}$ is non-degenerate,
  the usual polarization procedure produces a canonical compatible
  triple there.
 
  Recall that we assumed $r=\det(\omega)$ on an neighborhood $U$ $Z$,
  where $r$ is the coordinate function from Equation
  (\ref{eq:darboux_coords}). Since $\omega$ is folded symplectic,
  $\omega_E=r\cdot \mu$ for some non-degenerate positive 2-form $\mu$
  on $E$.  Therefore $A_E=r\cdot \tilde A_E$, where $\tilde
  A_E=g_E^{-1}\mu$ is a non-degenerate skew-endomorphism of $E$. By
  polarization of $\tilde A_E$ we obtain a smooth compatible triple
  $(\mu, h, \tilde J)$ on $E$, and using $A_E$ instead we get a
  compatible triple $(\omega_E, {g_J}_E, J_E)$ on $E\setminus E_Z$
  satisfying Equation (\ref{eq:triple_on_E}).
  
  On $Z$ we see that $\mu|_E=dr\wedge\alpha$ by Equation
  (\ref{eq:darboux}) and thus $\tilde A_E=R\otimes dr-\del_r\otimes
  \alpha=\tilde J$.
\end{proof}

The complex structure $\tilde J$ allows us to define a complement $K$
of $L$ in $E$ by $K=\tilde JL$, so we can refine the splitting
(\ref{eq:E-F-splitting}) over $Z$ to
\begin{eqnarray}\label{eq:K-L-F-splitting}
  T_ZX=K\oplus L\oplus F.
\end{eqnarray}

Equations (\ref{eq:triple_on_E}) show that $J$ is discontinuous across
the fold in the $E$ directions. However, on $U$ we may define two
smooth complex structures, denoted by $J^\pm$, such that
$J^\pm|_{X_\pm}=J$ by choosing $J^\pm_E=\pm \tilde J$.

We define folded holomorphic maps for all stable folded symplectic
4-manifolds, but at this stage most proofs rely on the following
additional structure that we assume along the fold.
\begin{definition}\label{def:S^1-invariant_fold}
  A stable folded symplectic manifold has a {\em circle invariant
    fold} (or $S^1$-invariant fold) if the flow of the associated
  characteristic vector field defines a free $S^1$ action on $Z$.
\end{definition}
Folded symplectic manifolds with $S^1$-invariant folds are especially
easy to work with.  They also occur frequently.  The standard folded
symplectic structure on the spheres (described in the next section) is
of this type. Connected sums of symplectic 4-manifolds always have
folded symplectic structures of this type (see \cite{ana1}). More
generally, connected sum along symplectomorphic symplectic
submanifolds of arbitrary codimension with symplectomorphic normal
bundles in symplectic 4-manifolds produces folded symplectic
structures with $S^1$-invariant folds. This can be seen by carefully
mimicking the symplectic connect sum construction in
\cite{mccarthy_wolfson} or \cite{IP_sum}.

$S^1$-invariant folds have a special structure, as described in Lemma
\ref{lem:pi_V-holomorphic} below. We will use this lemma repeatedly in
later sections. Whenever we work with $S^1$-invariant folds we assume
that the background metric $g'$ was chosen to be $S^1$-invariant on
the fold, resulting in an $S^1$-invariant compatible triple there.

\begin{lemma}\label{lem:pi_V-holomorphic}
  In the case of an $S^1$ invariant fold, $Z$ is an $S^1=\setR/\setZ$
  bundle over a symplectic manifold $(V,\omega_V)$ with projection
  $\pi_V:Z\rightarrow V$ such that
  \begin{enumerate}
  \item $\omega_Z=\pi_V^\ast\, \omega_V$
  \item there exists an $\omega_V$ compatible almost complex structure
    $J_V$ on $V$ such that $d\pi_V|_F$ is $(J,J_V)$ linear.
  \end{enumerate}
  Moreover, the 1-form $\alpha$ may be chosen such that
  \begin{eqnarray}\label{eq:curvature of alpha}
    d\alpha=C\cdot\omega_Z\qquad C= c_1(Z)/vol(V)
  \end{eqnarray}
  where $c_1$ is the first Chern class of the circle bundle
  $Z\rightarrow V$ and $vol(V)$ is taken with respect to $\omega_V$.
\end{lemma}

\begin{proof}
  Since the $S^1$ action on $Z$ is free we can exhibit $Z$ as an $S^1$
  bundle
  \begin{eqnarray*}
    \pi_V:Z\rightarrow V
  \end{eqnarray*}
  over a closed $(2n-2)$-dimensional manifold $V$. The kernel of
  $\omega$ coincides with the vertical subspace and
  $\mathcal{L}_R\omega=0$, so $\omega$ is $S^1$ invariant. Thus there
  exists a 2-form $\omega_V$ on $V$ such that $\omega=\pi_V^\ast
  \omega_V$. One readily checks that $\omega_V$ is non-degenerate and
  closed.
  
  Because $\ker(d\pi_V)$ is transverse to $F$,
  \begin{eqnarray*}
    d\pi_V|_{F}(z):F\rightarrow T_{\pi_V(z)}S
  \end{eqnarray*}
  is an isomorphism for each $z\in Z$. Since the complex structure $J$
  on $F$ is invariant under the $S^1$-action this map induces a
  complex structure $J_V$ on $V$ so that $d\pi_V|_{F}$ is $(J,J_V)$
  linear.
  
  To see equation (\ref{eq:curvature of alpha}) let $\alpha_0$ be a
  connection 1-form on $Z$, i.e. $\alpha_0$ is invariant under the
  $S^1$-action and satisfies $\alpha_0(R)=1$. Then $\iota_R
  d\alpha_0=\mathcal{L}_R\alpha_0=0$, so $d\alpha_0=\pi_V^\ast
  \omega_0$ is the pullback of a 2-from $\omega_0$ on $V$ which is
  just the curvature of the connection $\alpha_0$ and therefore
  represents $c_1(Z)$.
  
  Since $\omega_V$ is a volume form on $V$ there exists a constant
  $c\in\setR$ with
  \begin{eqnarray*}
    \int_V(c\,\omega_V-\omega_0)=0,
  \end{eqnarray*}
  so $(c\,\omega_V-\tilde\omega)=d\beta$ for some 1-form $\beta$ on
  $V$.  With the gauge transformation
  $\alpha=\alpha_0+\pi_V^\ast\beta$ we still have $\alpha(R)=1$, and
  \begin{eqnarray*}
    d\alpha=d\alpha_0+\pi_V^\ast
    d\beta=\pi_V^\ast\omega_0+\pi_V^\ast(c\omega_V-\omega_0)
    =c\,\pi_V^\ast\omega_V=c\,\omega_Z.
  \end{eqnarray*}
\end{proof}
Henceforth we assume that all folded symplectic manifolds in question
are compact, connected, 4-dimensional and have stable folds. If we
additionally assume the folds to be circle-invariant, we assume that
the 1-form $\alpha$ is chosen to satisfy Equation (\ref{eq:curvature
  of alpha}).

\vskip1cm
\section{Motivating Example}\label{sec:motivating_example}
To motivate what follows we will investigate possibilities to define
pseudo-holomorphic maps into $S^4$ with canonical folded symplectic
structure as defined below. One overruling principle is that we want
to obtain well-behaved moduli spaces of such maps. More precisely we
are looking for a notion of pseudo-holomorphic maps into folded
symplectic manifolds such that the linearized equations at a solution
give rise to a Fredholm operator and that the solutions are stable
under perturbations of the structures involved.

Recall that $S^4$ does not admit any symplectic form nor a continuous
almost complex structure. Therefore the answer to the question how to
generalize pseudo-holomorphic maps to this setting is far from
obvious. We let ourselves be guided by the folded symplectic structure
$\omega$. Since $\omega$ is non-degenerate on $X\setminus Z$, the
usual procedure to construct a compatible triple will yield an almost
complex structure $J$ there. Then it is clear what a
pseudo-holomorphic map from a Riemann surface $(\Sigma,j)$ into
$X\setminus Z$ is, namely a map with $(j,J)$-linear differential.
Since the fold $Z$ separates $X$ into $X_\pm$ this means that maps
from a connected domain into $X\setminus Z$ will have image in only
one side $X_+$ or $X_-$ of the fold. The question then is how to allow
maps to ``cross the fold'' i.e. have image on both sides of the fold.

One way to do this is to choose an almost complex structure on
$X\setminus Z$ that degenerates along $Z$ in a way that $X_\pm$ has
``cylindrical ends'' in the sense of Symplectic Field Theory
\cite{SFT}. This then reduces to the problem of holomorphic curves
relative to closed characteristics as discussed in \cite{IP_rel},
\cite{SFT} and \cite{SFT_compactness}.  In effect this is treating the
two sides $X_\pm$ as separate manifolds with boundary.

We try and find the analogue of holomorphic curves with non-degenerate
almost complex structure for folded symplectic manifold, while being
guided by the hope that these will limit to the relative curves
discussed above as we degenerate the almost complex structure. To do
this we define a circle-invariant folded symplectic structure on $S^4$
together with a compatible almost complex structure as in Lemma
\ref{lem:triple}.

View $S^4$ as the unit sphere in $\setR^5$. Then we have
\begin{itemize}
\item the restriction of the coordinate projection
  $\Pi:\setR^5=\setR\times\setR^4\rightarrow \setR^4$ to $S^4$, and
\item the stereographic projections $\sigma_\pm:S^4\setminus (\pm
  1,0,0,0,0)\rightarrow \setR^4$.
\end{itemize}

\begin{figure}[htbp]
  \centering
  \includegraphics[width=7cm]{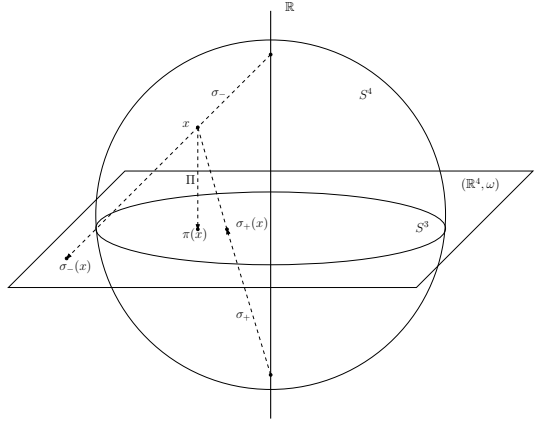}
  \caption{Folded Structures on $S^4$}
  \label{fig:folded S^4}
\end{figure}

Let $\omega_0$ and $J_0$ be the standard symplectic and complex
structures on $\setR^4$. Then $\omega=\Pi^\ast(\frac1\pi\omega_0)$ is
a folded symplectic form. The orientation induced by the folded
symplectic form agrees (disagrees) with the canonical orientation on
$S^4$ on the upper (lower) hemisphere $S^4_+$ ($S^4_-$), the fold
$Z=S^3$ is the intersection of $S^4$ with the equatorial plane
$\{(x_0,x_1,x_2,x_3,x_4)|x_0=0\}$. Choose the 1-form $\alpha$ on
$Z=S^3$ to be the canonical contact structure on $S^3$, i.e.  $\alpha$
is the restriction of the canonical 1-form
\begin{eqnarray*}
\alpha=\frac{1}{2\pi}\left(x_1dx_2-x_2dx_1+x_3dx_4-x_4dx_3\right)  
\end{eqnarray*}
to $S^3$. Thus $d\alpha=\omega_Z$.

With the almost complex structure $J$ defined as
$J|_{S^4_\pm}=\sigma_\pm^\ast J_0$ and the symmetric bilinear form
$g(u,v)=\omega(u,Jv)$ we obtain a folded compatible triple
$(J,\omega,g)$ on $S^4$. Note that all structures are invariant under
the $S^1=\setR/\setZ$-action on $Z$ given by the characteristic flow,
i.e. $t\cdot(z,w)=(e^{2\pi i\,t}z,e^{2\pi i\, t} w)$.

The map
\begin{eqnarray}\label{eq:involution}
  \tau:S^4\rightarrow S^4,\qquad
  (x_0,x_1,x_2,x_3,x_4)\mapsto(-x_0,x_1,x_2,x_3,x_4)
\end{eqnarray}
is a biholomorphic involution on $S^4$ exchanging the upper and lower
hemisphere and fixing the fold.

The bundle $F=\ker(\alpha)$ is given by the contact planes of the fold
$S^3$, and $E$ is spanned by the characteristic direction given by the
vertical subspaces of the Hopf fibration and the ``radial'' direction.

Any non-trivial $J$-holomorphic map from a Riemann surface
$(\Sigma,j)$ has to cross the fold, since each side is biholomorphic
to $B^4\subset\setC^2$ which is contractible and therefore does not
admit non-trivial holomorphic curves by their energy minimizing
property. Suppose for simplicity that the preimage of the fold is a
closed 1-dimensional separating submanifold $\sigma$ on the domain
$\Sigma$, cutting it into two parts $\Sigma_+$ mapping into $S^4_+$
and $\Sigma_-$ mapping into $S^4_-$.

Investigating the orientations on the transversal bundle $E$ near the
fold reveals that it is impossible to have holomorphic maps of class
$C^1$ (or even $C^0$) into $S^4$ that cross the fold (see also Remark
\ref{rem:orientations}). To remedy this one may look at curves that
are holomorphic on $\Sigma_+$, anti-holomorphic on $\Sigma_-$ and
continuous across the fold. Here we face the problem that too many
curves exists. For example, any holomorphic map from the upper
hemisphere of $S^2$ into $S^4_+$ with boundary on the fold can be
completed by an anti-holomorphic map from the lower hemisphere into
$S^4_-$ simply by reflection through the equator on domain and target,
i.e. using the anti-holomorphic involution of the hemispheres of $S^2$
and the holomorphic involution $\tau$ on $S^4$. Thus there exist an
infinite dimensional space of solutions in this case.

The reason for this is revealed by studying this as a boundary value
problem. The condition of ``continuous images'' does not give rise to
elliptic boundary conditions and therefore cannot lead to a Fredholm
problem. 

One possible remedy is to impose additional constraints to cut down
the solutions to a finite dimensional space. But there is no evident
way of doing this so that the structure of the solution space is
stable under perturbations and that allows for sufficient interaction
between the two sides of the fold.

Another approach is to allow discontinuous images. We think of this as
holomorphic curves that leave the fold at a location that is different
from where they enter, the relation between these is given by a
tunneling map in the fold. This way we define a Fredholm problem for
discontinuous pseudo-holomorphic maps into folded symplectic
manifolds. Roughly speaking, a folded holomorphic map consists of
\begin{itemize}
\item a domain $(\Sigma,j)$ with 1-dimensional submanifold $\sigma$
  separating $\Sigma$ into $\Sigma_+$ and $\Sigma_-$
\item $J$-holomorphic maps $u_+:(\Sigma_+,\sigma)\rightarrow (X_+,Z)$
  and $u_-:(\Sigma_-,\sigma)\rightarrow (X_-,Z)$
\item and a pair of tunneling maps in $Z$ connecting $u_+(\sigma)$ to
  a closed characteristic and then continuing on to connect to
  $u_-(\sigma)$.
\end{itemize}

\begin{figure}[htbp]
  \centering
  \includegraphics[width=7cm]{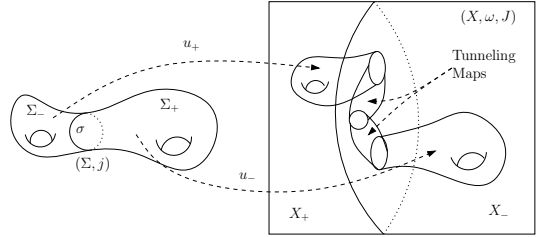}
  \caption{The map tunnels through the fold, exiting the fold at a
    location that is different from where it entered.}
  \label{fig:folded maps 1}
\end{figure}

We will make this precise in the following sections. In order to
convince the reader that it is natural to consider discontinuous maps
in the context of folded symplectic manifolds we will first give a
trivial example.

Consider the complex elliptic fibration $X=E(1)$ and let
$T\hookrightarrow E(1)$ be a regular fiber. Perturb the almost complex
structure $J$ so that a neighborhood $X_-$ of $T$ is biholomorphic to
$X_-=D\times T$ with a product complex structure, where $D$ is the
closed unit disk in $\setC$ with the canonical complex structure and
$T=\{0\}\times T$ with its original complex structure. Set
$X_+=E(1)\setminus X_-$.

After choosing $z_0\in S^1$ we define for $z\in S^1$ the element
$\overline z\in S^1$ by reflection through $z_0$, i.e.
$z\cdot\overline z=z_0$. Using the two self-maps of the boundary
$Z=\del X_\pm=S^1\times T$
\begin{eqnarray*}
  \iota(z,w)=(z,w)\qquad\mathrm{and}\quad
  \iota_F(z,w)=(\overline z, w)
\end{eqnarray*}
we construct the manifolds
\begin{eqnarray*}
  E(1)=X_+\sqcup_{\iota} X_-\qquad\mathrm{and}\quad
  E^F(1)=X_+\sqcup_{\iota_F} \overline{X_-}
\end{eqnarray*}
with almost complex structure inherited from each piece.  With a
little more care we can ensure that $E^F(1)$ comes equipped with a
folded symplectic structure, independent of the choice of $z_0$.

Consider the biholomorphism
\begin{eqnarray*}
  \Psi:E(1)\setminus Z\rightarrow E^F(1)\setminus Z
\end{eqnarray*}
given by the identity map on each piece, $X_+$ and $X_-$.

This setup suggest the following definition for folded holomorphic
maps into the folded symplectic manifolds $E^F(1)$. Let $\sigma$ be a
separating submanifold in a Riemann surface $\Sigma$, separating
$\Sigma$ into two pieces $\Sigma_+$ and $\Sigma_-$.
\begin{definition}\label{def:folded_hol_E(1)}
  A folded holomorphic map consists of a pair of maps $(u_+,u_-)$ such
  that 
  \begin{eqnarray*}
    u_\pm:\Sigma_\pm\rightarrow X_\pm,\qquad 
    u_\pm=\Psi\circ u|_{\Sigma_\pm}
  \end{eqnarray*}
 for some $J$-holomorphic $u:\Sigma\rightarrow
  E(1)$.
\end{definition}
Note that there is an $S^1$ ambiguity in the gluing map $\iota_F$, so
there is really an $S^1$ family of these folded holomorphic maps.  It
is clear from the definition that this will yield a well-defined
moduli space, although the maps are necessarily discontinuous.

\begin{figure}[htbp]
  \centering
  \includegraphics[width=9cm]{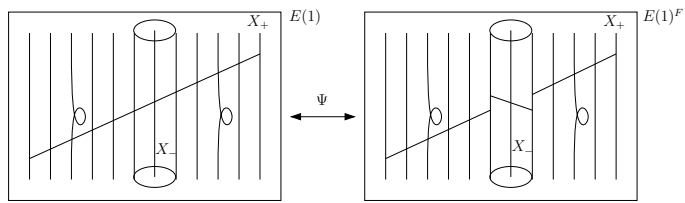}
  \caption{Folded holomorphic maps in $E^F(1)$.}
  \label{fig:E(1)_2}
\end{figure}

In the following we will show how to define folded holomorphic maps
into stable folded symplectic manifolds that reduce to Definition
\ref{def:folded_hol_E(1)} in the case of $E^F(1)$.

Note that the pieces $u_+$ and $u_-$ of a folded holomorphic map into
$E^F(1)$ have boundary values in
\begin{eqnarray*}
  \Delta^Z=\{( u_+|_\sigma, u_-|_\sigma)\}=\graph(\Psi)\in 
  \Map(\sigma,Z)\times \Map(\sigma,Z)
\end{eqnarray*}
We call $\Delta^Z$ the {\em folded diagonal}. The above definitions
work in the case where the fold $Z$ has the structure of a trivial
$S^1$-bundle with the characteristic foliation being vertical. To
define folded holomorphic maps into more general folded symplectic
manifolds we need to generalize the folded diagonal so that it
continues to give elliptic boundary conditions for the maps $u_+$ and
$u_-$.

\vskip1cm
\section{Folded Holomorphic Maps}\label{sec:folded maps}
We define folded holomorphic maps and lay the functional analytic
foundation for the later sections. We start by describing the domains
of folded holomorphic maps, then we set up the Sobolev spaces and
lastly we set up the PDE. 

\subsection{Folded Domains}
\label{sec:domains}
We define the domains of folded holomorphic maps. Due to the
additional structure needed for the tunneling maps this is more
involved than one might at first expect.

\begin{definition}[Folded Domain]\label{def:folded domain}
  A {\em smooth rational folded domain} $\mathcal{D}$ consists of
  \begin{enumerate}[(i)]
  \item closed Riemann surfaces $(\Sigma_0, j_0, {\bf p}_0)$ and
    $(\Sigma_1, j_1, \{p\})$ of genus $g_0$ and $g_1$ with
    ${\bf p}_0=\{p_1,\ldots p_n\}\subset\Sigma_0$
  \item continuous functions $\tau_i:\Sigma_i\rightarrow\setR$ so that
    with
    \begin{eqnarray*}
      \Sigma_i^+=\tau_i^{-1}[0,\infty)\qquad
      \Sigma_i^-=\tau_i^{-1}(-\infty,0]
    \end{eqnarray*}
    their restrictions $\tau_i^\pm=\tau_i|_{\Sigma_i^\pm}$ extend to
    smooth functions on $\Sigma_i$ with zeros of at most finite order
    along $\sigma_i=\tau_i^{-1}(0)=(\tau_i^\pm)^{-1}(0)$ with matching
    order of vanishing at each point of $\sigma_i$, $i=1,2$.
  \item $ p\in\Sigma_1^-$
  \item a function $g:\Sigma_0^+\rightarrow \setR$ and a
    diffeomorphism $\psi:\Sigma_0^+\rightarrow\Sigma_1^+$ satisfying
    the conditions
    \begin{eqnarray}\label{eq:folded domain condition}
      \psi^\ast j_1=j_0,\qquad \psi^\ast \tau_1=e^g\tau_0.
    \end{eqnarray}
  \end{enumerate}
\end{definition}
We set ${\bf p}_0^\pm={\bf p}_0\cap\Sigma_0^\pm$ to be the marked
points contained in $\Sigma_0^\pm$. The zero sets
$\sigma_0$ and $\sigma_1$ are called the
{\em domain folds}. The map $\psi$ gives an identification of
$\sigma_0$ and $\sigma_1$, so we will refer to the domain folds simply
by $\sigma$. To simplify notation we set
\begin{eqnarray*}
\Sigma_\pm=\Sigma_0^\pm,\qquad S=\Sigma_1^-,\qquad \dot S=S\setminus
\{p\}.
\end{eqnarray*}
When no confusion can occur we will drop the subscripts on $j_i$ and
$\tau_i$.

\begin{figure}[htbp]
  \centering
  \includegraphics[width=9cm]{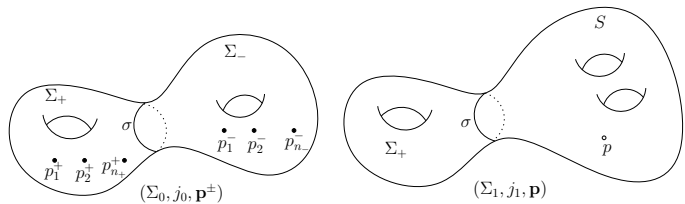}
  \caption{Folded Domains}
  \label{fig:folded domains}
\end{figure}

The purpose of the functions $\tau_0$ and $\tau_1$ is to give possible
locations of the domain fold $\sigma$, but we are not interested at
the specific values of $\tau_i$ away from there zero set.

\begin{definition}\label{def:gauge_group}
  The group
  \begin{eqnarray}
    \label{eq:gauge group}
    \mathcal{G}=\Map(\Sigma_1,\setR)
    \times\Diff^+(\Sigma_0,{\mathbf p}_0)
    \times\Diff^+(\Sigma_1,\{p\})
  \end{eqnarray}
  is called the {\em gauge group}.
\end{definition}

\begin{lemma}\label{lem:domain_gauge_action}
  The gauge group $\mathcal{G}$ acts on the space of folded domains in
  the following way. Given $(f,\phi_0,\phi_1)\in\mathcal{G}$ the
  folded domain transforms according to
  \begin{eqnarray*}
    j_i&\mapsto& \phi_i^\ast j_i=d\phi_i^{-1}\circ j_i\circ d\phi_i,\\
    \tau_0&\mapsto& \tau_0\circ \phi_0,\qquad\qquad
    \tau_1\mapsto e^f\tau_1\circ \phi_1,\\
    {\bf p}_0&\mapsto& \phi_0^\ast({\bf p}_0),\qquad\qquad
    p\mapsto \phi_1^\ast(p),\\
    \psi&\mapsto&\phi_1^{-1}\circ\psi\circ \phi_0,\\
    g&\mapsto&g\circ \phi_0+(\phi_1^{-1}\circ\psi\circ \phi_0)^\ast f.
  \end{eqnarray*}
\end{lemma}

\begin{remark}
{\em
    \begin{enumerate}
    \item The sets ${\bf p}_0^\pm$ are disjoint unless some
      of the marked points lie on $\sigma_0$.
    \item Let be $\hat S$ the radial compactification of $\dot S$ as
      in \cite{SFT}. Then we define the associated 2-dimensional
      $C^0$-cycles with boundary $\tilde \Sigma_\pm$ by gluing
      $\Sigma_+$ to $\hat S$ to get $\tilde \Sigma_+=\hat\Sigma_1$ or
      by gluing $\Sigma_-$ to $-\hat S$ to get $\tilde \Sigma_-$,
      where the gluing is done along $\sigma$.  $\tilde \Sigma_+$ is a
      smooth Riemann surface, whereas the same only holds for
      $\tilde\Sigma_-$ in the case that the domain fold $\sigma$ is a
      manifold.
    \item Since the order of vanishing of
      $\tau^\pm=\tau_0|_{\Sigma_\pm}$ at each point of $\sigma$ is the
      same we may define the ``gap-function''
      \begin{eqnarray*}
        a:\sigma\rightarrow (0,\infty)
      \end{eqnarray*}
      at a point $z\in\sigma$ by choosing a path
      $\gamma:[-1,1]\rightarrow \Sigma$ with
      $\gamma([0,1])\in\Sigma_+$ and $\gamma([-1,0])\in\Sigma_-$ and
      $\gamma(0)=z$ and setting
      \begin{eqnarray}
        a(z)=\lim_{t\rightarrow 0}
        \frac{\tau_-\circ\gamma(|t|)}{\tau_+\circ\gamma(-|t|)}.
        \label{eq:a_tau}
      \end{eqnarray}
      Note that this is independent of the choice of path $\gamma$.
      At a generic point $z\in\sigma$ where $\tau_\pm(z)$ vanish to
      first order we have
      \begin{eqnarray}\label{eq:a_tau2}
        a(z)=\frac{d\tau_-(z)}{d\tau_+(z)},
      \end{eqnarray}
      where the quotient of forms is to be understood as being taken
      after evaluating on an arbitrary normal vector to $\sigma$. This
      is well defined since $d\tau_\pm$ are collinear at $z$ as they
      vanish along the direction of $\sigma$.
  \end{enumerate}
}
\end{remark}

\subsection{Sobolev Spaces of Maps}
\label{sec:sobolev}
To proceed we need to give a precise definition of the Sobolev spaces
we plan to use. We work with the non-degenerate metric $\tilde g$
that, with the notation of Lemma \ref{lem:triple}, leaves the
splitting $TU=E\oplus F$ invariant, equals $g$ on $F$ and $\mu$ on $E$
in a neighborhood $U'\subset \bar U'\subset U$ of the fold and equals
$g$ outside of $U$ so that $J_\pm$ is antisymmetric w.r.t. $\tilde g$.
In particular $\tilde g=dr\otimes dr+\alpha\otimes \alpha+g_F$ on $Z$.

Consider a folded domain $\mathcal D$ as in Definition \ref{def:folded
  domain}. We follow the definitions from \cite{booss_wojchiechowski}
for Sobolev spaces on manifolds with boundary.
\begin{definition}
  Fix a Riemannian metric in the conformal class of $j$ on $\Sigma_0$
  and positive integers $k,p$ with $kp>2$. Let $U_\pm\subset\Sigma_0$
  be open subsets properly containing $\Sigma_\pm$ and that are of the
  same homotopy type as $\Sigma_\pm$. Then let
  $W^{k,p}(\Sigma_\pm,X_\pm)$ be the Banach manifold consisting of
  maps $f:\Sigma_\pm\rightarrow X_\pm$ that are restrictions of maps
  $\tilde f:U_\pm\rightarrow X$ of class $W^{k,p}$ to $\Sigma_\pm$
  that send the domain fold $\sigma$ into $Z$.
\end{definition}

$W^{k,p}(\Sigma_\pm,X_\pm)$ is a smooth separable Banach manifold
modeled locally at a map $u_\pm\in W^{k,p}(\Sigma_\pm,X_\pm)$ on the
space $W^{k,p}(u_\pm^\ast TX)$.

Next we define the Banach manifolds of maps from punctured surfaces
into $Z$. Here we differ from the traditional treatment found in
\cite{schwarz}, \cite{hofer_fredholm} or \cite{dragnev} in the basic
definitions as we do not a priory specify the asymptotics at the
punctures but rather allow the maps to converge to arbitrary closed
characteristics. We also avoid using the auxiliary $\setR$-factor in
the ``symplectization'' as it does not add any information. We believe
that this approach gives a more natural setup for the Fredholm theory
needed in our case.

For a closed Riemann surface $\Sigma$ with finitely many punctures
$\{p_k\}$ and $\dot \Sigma=\Sigma\setminus \{p_k\}$ we let
$W^{k,p}_{loc}(\dot\Sigma,Z)$ be the space of maps from $\dot\Sigma$
to $Z$ that, in local coordinates, are in
$W^{k,p}_{loc}(\setR^2,\setR^3)$.

For a Riemann surface $\Sigma$ with boundary and finitely many
punctures we assume that $\Sigma\subset \Sigma'$ for some open
Riemannian manifold $\Sigma'$ of the same homotopy type and we set
$W^{k,p}_{loc}(\dot\Sigma,Z)$ to be the space of maps from
$\dot\Sigma$ to $Z$ that are restrictions of maps in
$W^{k,p}_{loc}(\dot\Sigma',Z)$.

Let $(\Sigma,j,\{p_k\})$ be a Riemann surface (possibly with boundary)
with conformal structure $j$ and punctures $\{p_k\}$. Set
$\dot\Sigma=\Sigma\setminus\{p_k\}$.  For $r\in \setR$ define the
half-infinite cylinder
\begin{eqnarray*}
  C_r=[r,\infty)\times S^1
\end{eqnarray*}
with coordinates $s\in[r,\infty)$, $t\in S^1=\setR/\setZ$ complex
structure $j$ with $j\del_s=\del_t$ and volume form $\dvol=ds\wedge
dt$.  Set $C=C_0$. Then at each puncture $p_k$ we have local conformal
coordinates $\sigma_k:C\rightarrow \dot\Sigma$.

Fix a constant $\delta> 0$. Since we will be only interested in
$\delta$ close to zero we assume throughout that $\delta$ is bounded
by some constant $M$.
\begin{definition}[Asymptotic Energy]
  For maps $v\in W^{k,p}_{loc}(C,Z)$ we define the {\em asymptotic
    energy}
  \begin{eqnarray}\label{eq:asymptotic_energy}
    E_r(v)=\int_{C_r} (|v^\ast\alpha(\del_s)|^2
    +|d(v^\ast\alpha(\del_t))|^2
    +|\pi_F\, dv|^2)\;e^{\delta s}\dvol.
  \end{eqnarray}
\end{definition}

\begin{definition}\label{def:tunneling_sobolev}
  Let $W^{k,p}_\delta(C,Z)$ be the space of finite asymptotic energy
  $W^{k,p}_{loc}(C,Z)$ maps. Similarly, let
  $W^{k,p}_\delta(\dot\Sigma,Z)$ be the space of
  $W^{k,p}_{loc}(\dot\Sigma,Z)$ maps that are in $W^{k,p}_\delta(C,Z)$
  in some local conformal coordinates at each puncture.
\end{definition}

\begin{definition}
  For $v\in W^{k,p}_\delta(C,Z)$ let $W^{k,p}_\delta(C,v^\ast TZ)$ to
  be the space of sections $\zeta\in W^{k,p}_{loc}(C,v^\ast TZ)$ that
  have finite asymptotic energy, i.e. that satisfy
  \begin{eqnarray}\label{eq:variational_energy}
    E_r(\zeta)=\int_{C_r} \left(|\pi_F(\nabla\zeta)|^2
      +|\alpha({\nabla}_s\zeta)|^2
      +|d(\alpha({\nabla}_t\zeta))|^2\right)e^{\delta s}\dvol
    <\infty.
  \end{eqnarray}
  Analogously, for $v\in W^{k,p}_\delta(\dot \Sigma,Z)$ we define
  $W^{k,p}_\delta(\dot \Sigma,v^\ast TZ)$ to be the space of sections
  $\zeta\in W^{k,p}_{loc}(\dot\Sigma,v^\ast TZ)$ that are in
  $W^{k,p}_\delta(C,v^\ast TZ)$ in local conformal coordinates at each
  puncture.
\end{definition}

With these definitions $W^{k,p}_\delta(\dot\Sigma,Z)$ is
a separable Banach manifold, modeled locally at a map $v\in
W^{k,p}_\delta(\dot\Sigma,Z)$ on a neighborhood of the zero
section in $W^{k,p}_\delta(\dot\Sigma,v^\ast TZ)$.

\subsection{Folded Maps and Folded Holomorphic
  Maps}\label{sec:folded_maps_2}

\begin{definition}[Space of Folded Maps]\label{def:folded_maps}
  Fix positive integer $k$ and $p\in\setR$ with $kp>2$, non-negative
  integers $g_0,\,g_1,\,n$ and relative homology classes $A_\pm\in
  H_2(X_\pm,\mathcal{R};\setZ)$. Then a folded map $(u_+,u_-)$ with
  respect to $A_\pm$ consists of
  \begin{enumerate}[(i)]
  \item a folded domain and
  \item $u_\pm\in W^{k,p}(\Sigma_\pm,X_\pm)$
  \end{enumerate}
  such that there exist maps $v_\pm\in W_\delta^{k,p}(\dot S,Z)$ with
  \begin{eqnarray}\label{eq:pullback_det}
    u_\pm^\ast \det(\omega)=\tau,\qquad 
    u_\pm|_\sigma=v_\pm|_\sigma,\qquad\mathrm{and}\quad
    \left[u_\pm\sqcup_\sigma v_\pm\right]=A_\pm.
  \end{eqnarray}
\end{definition}

\begin{figure}[htbp]
  \centering
  \includegraphics[width=7cm]{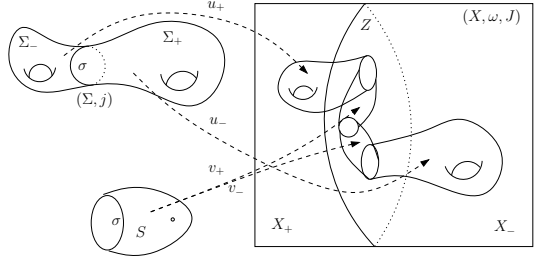}
  \caption{Folded (Holomorphic) Maps}
  \label{fig:folded holomorphic maps}
\end{figure}

\begin{definition}[Folded Holomorphic
  Maps]\label{def:folded_holomorphic_map}
  A folded map $(u_+,u_-)$ is called {\em folded holomorphic} if
  $\delbar_J u_\pm=0$ and $(u_+,u_-)|_\sigma\in \Delta^Z$, where the
  folded diagonal $\Delta^Z$ is defined in \ref{def:folded_diagonal}.
\end{definition}

The key to defining folded holomorphic maps lies in the
definition of the ``matching condition''. It is specified by the
subspace
\begin{eqnarray*}
  \Delta^Z\subset\Map(\sigma,Z)\times \Map(\sigma,Z),
\end{eqnarray*}
called the {\em folded diagonal}. Intuitively, we view $\Delta^Z$ as a
scattering function that takes boundary conditions of ``incoming''
holomorphic curves on the ``$+$'' side and transforms them into
boundary conditions for ``outgoing'' holomorphic curves on the ``$-$''
side. To make sense out of this we don't need to define this
scattering function for every element of $\Map(\sigma,Z)$ but it
suffices to define it on the maps that are possible boundary
conditions of holomorphic maps. The precise definition of the folded
diagonal $\Delta^Z$ is rather involved and we postpone it until
Section \ref{sec:def_tunneling}.

\begin{lemma}\label{lem:gauge_action}
  The gauge group $\mathcal{G}$ defined in Definition
  (\ref{def:gauge_group}) acts on the space of folded holomorphic maps
  by precomposition of $u_\pm$ by the element in $\Diff^+(\Sigma_0)$,
  on $v_\pm$ by the element in $\Diff^+(\Sigma_1)$ and acts on the
  folded domain as in Lemma \ref{lem:domain_gauge_action}.
\end{lemma}

\begin{proof}
  The action preserves solutions to the holomorphic map equation and
  the condition that $\tau_0=u_\pm^\ast\det(\omega)$, since $\tau_0$
  and $u_\pm$ are acted on by precomposition with the same
  diffeomorphism. Moreover it preserves the identification of the
  domain folds $\sigma_0$ and $\sigma_1$ as for $p\in\sigma$
  \begin{eqnarray*}
    u_\pm'(p)
    =u_\pm\circ\phi_0(p)
    =v_\pm\circ\psi\circ\phi_0(p)
    =v_\pm'\circ\phi_1^{-1}\circ\psi\circ\phi_0(p)
    =v_\pm'\circ\psi'(p).
  \end{eqnarray*}
  So this action preserves the set of folded holomorphic maps as long
  as it preserves the folded diagonal. We will postpone this part of
  the proof until Remark \ref{rem:gauge_invariant}.
\end{proof}

We make an important observation about the orientations of maps along
the fold.
\begin{remark}\label{rem:orientations}
  {\em Let $p\in\sigma$ be a point where $\tau$ vanishes transversely
    and let $\eta\in T_p\Sigma$ be an outward normal vector to
    $\sigma$ at $p$, i.e. $j\eta\in T_p\sigma$ and
    $d\tau_\pm(\eta)<0$.  Recall the notation from Lemma
    \ref{lem:triple}. Then, remembering that
    $\tau_\pm=u_\pm^\ast\det(\omega)$ and $J_\pm\del_r=\pm R$ and thus
    $dr=\pm \alpha\circ J_\pm$,
    \begin{eqnarray*}
      u_\pm^\ast\alpha(j\eta)&=&\alpha(du_\pm j\eta)
      =\alpha\circ J_\pm (du_\pm\eta)
      =\pm dr(du_\pm(\eta))
      =\pm d(\det(\omega))(du_\pm(\eta))\\
      &=&\pm d\tau_\pm(\eta)
    \end{eqnarray*}
    so the values of $u_\pm^\ast\alpha$ on tangent vectors have
    opposite sign.
    
    This shows that at a point $p\in \sigma$ where $\tau_\pm$ vanish
    transversely, the gap-function $a$ from Equation (\ref{eq:a_tau})
    can be expressed in terms of $u_\pm$ using Equation
    (\ref{eq:a_tau2}) as
    \begin{eqnarray}
      a(p)=-\frac{u_-^\ast\alpha(p)}{u_+^\ast\alpha(p)}
      \label{eq:a_u}
    \end{eqnarray}
    where the fraction is to be understood as evaluated on an
    arbitrary element of $T_p\Sigma$ such that the denominator does
    not vanish. This is well defined since $u_\pm^\ast\alpha$ are
    collinear at $p$.
    }
\end{remark}
This turns out to be a crucial observation in showing that the
folded diagonal poses an elliptic boundary condition and that
therefore the linearized operator is Fredholm.

\vskip1cm
\section{Tunneling Maps}\label{sec:tunneling}
Throughout this section we assume that the function $\tau_\pm$
vanishes transversely, so the domain fold $\sigma$ is a manifold.

Tunneling maps give the matching conditions for pseudo-holomorphic
maps into $X_\pm$, add homological data to the maps, and ensure that
families of folded holomorphic maps have bounded energy. They are
central in proving regularity of solutions and they guarantee that the
linearized operator is Fredholm. In this section we define all
relevant structures and discuss their properties.

\subsection{Definition of Tunneling Maps and the Folded Diagonal}
\label{sec:def_tunneling}
Tunneling maps are maps into the fold $Z$ that connect the images of
the folded maps into $X_+$ and $X_-$. All tunneling maps will have
domains $(\dot S,j)$ of the form $(\Sigma_1^-\setminus\{p\},j_1)$ with
boundary $\del S=\sigma$.

Tunneling maps satisfy an equation that depends only on the CR
structure $(Z,F,J)$ on the fold $Z$ and the 1-form $\alpha$. First we
need the following generalization of the holomorphic map equation.
\begin{definition}[$\mathcal{H}$-Holomorphic Maps]\label{def:H_hol}
  A map $v:\dot S\rightarrow Z$ is called {\em
  $\mathcal{H}$-holomorphic} if
  \begin{eqnarray}
    \delbar^F_J v
    &=&\frac{1}{2}\left(\pi_F\,dv+J\,\pi_F\,dv\,j\right)=0
    \label{eq:H_F}\\
    d(v^\ast\alpha\circ j)&=&0.\label{eq:H_L}
  \end{eqnarray}
\end{definition}
$\mathcal{H}$-holomorphic maps can alternatively viewed as families of
$J$-holomorphic maps into the complex cylinder (``symplectization'')
over $Z$ parametrized by
\begin{eqnarray*}
  \mathcal{H}_n(S,j)=\{\lambda\in\Omega^1(S)|d\lambda=0,\
  d^\ast\lambda=0,\ \mathrm{and}\ \lambda(j\eta)=0\,\forall\eta\in
  T\del S\}.
\end{eqnarray*}
$\mathcal{H}_n$ is naturally isomorphic to $H^1(S,\setR)$ as a vector
space and consists of smooth 1-forms (see \cite{duff_spencer}).  In
\cite{planar_weinstein} the authors called maps from punctured
surfaces without boundary satisfying this condition ``generalized
holomorphic''.

\begin{definition}[Tunneling Maps]\label{def:tunneling}
  Fix a positive integer $k$ and real $p>0$ such that $kp>4$. A
  tunneling map is a $\mathcal{H}$-holomorphic map of class
  $W^{k,p}_\delta(\dot S,Z)$.
\end{definition}

We make a simple, but essential observation about tunneling maps.
\begin{lemma}\label{lem:no periods}
  Let $v$ be a tunneling map. Then the periods of $v^\ast\alpha\circ
  j$ vanish in a neighborhood of the puncture.
\end{lemma}
\begin{proof}
  Fix local conformal coordinates $C=[0,\infty)\times S^1$ at the
  puncture as in Definition \ref{def:tunneling_sobolev}.  Since
  $v^\ast\alpha\circ j$ is closed, the value of
  \begin{eqnarray*}
    \int_{\{r\}\times S^1}v^\ast\alpha\circ j
    =\int_{\{r\}\times S^1}v^\ast\alpha(\del_s)dt
  \end{eqnarray*}
  does not depend on $r$. Observing that $A_r=(r,r+1)\times S^1$ has
  unit area and using equation (\ref{eq:variational_energy}) we
  compute
  \begin{eqnarray*}
    \left|\int_{\{r\}\times S^1}v^\ast\alpha\circ j\right|
    \le\int_{A_r}|v^\ast\alpha(\del_s)|\dvol
    \le e^{-\delta r/2}\left(\int_{A_r}e^{\delta
        s}|v^\ast\alpha(\del_s)|^2\dvol\right)^{\frac12}
    \le e^{-\delta r/2}\sqrt{E_r(v)},
  \end{eqnarray*}
  showing that the periods of $v^\ast\alpha\circ j$ are arbitrarily
  small and therefore vanish.
\end{proof}

Tunneling maps are very well-behaved. They satisfy an elliptic system
of PDEs and are therefore smooth on the interior of $\dot S$. Given
elliptic boundary conditions their linearized operator is Fredholm for
almost all choices of weight $\delta$ by Theorem 6.3 of
\cite{lockhart_mcowen}.  They have nice limits at the punctures,
namely they converge to a closed characteristic exponentially fast.
Thus they extend to continuous maps from the radial
compactification of the domain. We will blur the distinction between a
tunneling map and its continuous extension to the radial
compactification.

\begin{remark}\label{rem:action of gauge group}
  {\em Note that the asymptotic energy (\ref{eq:asymptotic_energy}) is
    invariant under diffeomorphisms of the domain, so the gauge group
    $\mathcal{G}$ from Definition \ref{def:gauge_group} acts on the
    space of tunneling maps by precomposition.  }
\end{remark}

\subsection{Conjugate Tunneling Maps}
\label{sec:conjugate_tunneling}
Tunneling maps connect boundary values of $J$-holomorphic maps in
$X_\pm$ to closed characteristics. To get a scattering function of
incoming boundary values from a map $u_+$ into $X_+$ to outgoing
boundary values of a map $u_-$ into $X_-$, we need to define a
relation between tunneling maps $v_+$ connecting to $u_+$ and $v_-$
connecting to $u_-$.

We clarify some notation used in this section. For a tunneling map $v$
and an element $m\in S^1(T_p S)$ in boundary $\hat S\setminus \dot S$
of the radial compactification $\hat S$ of $\dot S$ we have that
$v(m)$ has image on some closed characteristic $x\in\mathcal{R}$.
Given a parametrization $x(t):S^1\rightarrow Z$ of the closed
characteristic, it inherits the group structure from $S^1=\setR/\setZ$
via $x(s)+x(t)=x(s+t)$.

\begin{definition}[Conjugate Tunneling Maps]
  \label{def:conjugate_tunneling}
  Two tunneling maps $v_+$ and $v_-$ with domain $(\dot S,j)$ that are
  asymptotic to the same closed characteristic $x\in\mathcal{R}$ are
  called {\em conjugate} with respect to a parametrization $x(t)$ of
  $x$ if there exists a function $f:S\rightarrow (0,\infty)$ with
  $f|_\sigma\equiv 1$ such that
  \begin{eqnarray}
    v_+^\ast\omega&=&f v_-^\ast\omega
    \qquad\mathrm{on}\ \dot S\label{eq:conjugate_F}\\
    \lambda=v_+^\ast\alpha\circ j+v_-^\ast\alpha\circ j&=&0 
    \qquad \mathrm{on}\ T\sigma\label{eq:conjugate_L}\\
    v_+(m)+v_-(m)&=&x(0)\qquad\forall\ m\in S^1(T_pS)\in\hat S.
    \label{eq:conjugate_marker}
  \end{eqnarray}
  
  We moreover demand that $v_+$ and $v_-$ converge to the same
  eigenvector of the asymptotic operator (c.f.
  \cite{planar_weinstein}, Equation (2.3)) at the puncture.
  
  Let
  \begin{eqnarray}
    \label{eq:Q}
    Q=\{q\in S|v_\pm\ \mathrm{fails\ to\ be\ an\ immersion\ at\ }q\}.
  \end{eqnarray}
  In the case of $S^1$-invariant folds we further demand that the
  order of vanishing of the 1-form $\lambda$ defined in
  (\ref{eq:conjugate_L}) is greater than or equal to the one of
  $\pi_F\,dv_\pm$ at each point $q\in Q$.
\end{definition}

\begin{remark}\label{rem:conjugate}
  {\em If $v_+$ and $v_-$ are conjugate, then their projections
    $\tilde v_\pm=\pi_V\circ v_\pm$ to the base $V$ agree by Equation
    (\ref{eq:conjugate_F}) and the requirement that $v_\pm$ converge
    to the same eigenvector of the asymptotic operator. To see this
    look at the projections $\tilde v_\pm=\pi_V\,v_\pm$. By Equation
    (\ref{eq:conjugate_F}) their differentials satisfy $d\tilde
    v_+=\tilde f\cdot\tilde dv_-$ on $S$ for some nowhere vanishing
    holomorphic function $\tilde f:S\rightarrow \setC$. Since the
    function $f$ in Equation (\ref{eq:conjugate_F}) equals 1 on
    $\sigma$ we conclude that $\tilde f|_\sigma:\sigma\rightarrow S^1$
    takes values in the unit circle. Since $\tilde f$ has degree 0 we
    conclude that $\tilde f$ is constant. The requirement that $v_\pm$
    converge to the same eigenvalue of the asymptotic operator implies
    that $\tilde f(p)=1$ and since $v_\pm$ converge to the same closed
    characteristic we have $\tilde v_+=\tilde v_-$.
    
    Thus there exists an $S^1$-valued function
    \begin{eqnarray*}
      g:\dot S\rightarrow S^1\qquad\mathrm{with}\ v_-(z)=g(z)\ast v_+(z)
    \end{eqnarray*}
    where $\ast$ denotes the $S^1$ action on $Z$ given by the
    characteristic flow. Equation (\ref{eq:conjugate_marker}) implies
    that $v_+$ and $v_-$ have opposite multiplicities at the punctures
    and the degree of $g$ at the puncture is $-2d$ if $d$ is the
    degree of $v_+$ at the puncture.
    
    $g$ satisfies $dg=v_-^\ast\alpha-v_+^\ast\alpha$ so 
    \begin{eqnarray*}
      d^\ast dg&=&\ast d(dg\circ j)=\ast d(v_-^\ast\alpha\circ
      j-v_+^\ast\alpha\circ j)=0,\\
      dg\circ j&=&-2v^\ast\alpha\circ j\qquad\mathrm{on}\ T\del S
    \end{eqnarray*}
    where we used Equations (\ref{eq:H_L}) and
    (\ref{eq:conjugate_L}). Therefore $g$ is a harmonic $S^1$-valued
    function with prescribed asymptotics satisfying von Neumann
    boundary conditions. This determines $g$ uniquely up to
    $H^1(S;\setZ)$ and an overall constant. Equation
    (\ref{eq:conjugate_marker}) this fixes the overall constant. By
    Remark \ref{rem:orientations} we conclude that only a finite
    number, bounded by a constant only depending on $d$, of values for
    the periods of $dg$ are possible.
    
    In the circle invariant case the zeros of $\pi_F\,dv_\pm$
    automatically agree. Thus the dimension of the space of conjugate
    tunneling maps jumps as the tunneling maps acquire points $q\in
    Q$. To deal with this non-generic situation we introduced the
    additional restriction on the order of vanishing of $\lambda$
    along $Q$. We will see that this ensures transversality for
    generic conjugate tunneling maps.
  }
\end{remark}

We adopt the convention that for folded maps $(u_+,u_-)$ and tunneling
maps $(v_+,v_-)$ their restrictions to $\sigma$ are denoted by a hat,
i.e.
\begin{eqnarray*}
  \hat u_\pm=u_\pm|_\sigma,\qquad\mathrm{and}\quad
  \hat v_\pm=v_\pm|_\sigma.
\end{eqnarray*}
We are now prepared for the main definition.
\begin{definition}[Folded Diagonal]\label{def:folded_diagonal}
  The {\em folded diagonal} is the subset in $\Map(\sigma,Z)\times
  \Map(\sigma,Z)$ defined by
  \begin{eqnarray*}
    \Delta^Z=\left\{(\hat v_+,\hat v_-)\big|
        (v_+,v_-)\ \mathrm{are\ conjugate} \right\}.
  \end{eqnarray*}
\end{definition}

\begin{remark}\label{rem:gauge_invariant}
  {\em Note that the folded diagonal is invariant under the action of
    the gauge group $\mathcal{G}$. Indeed, if $(v_+,v_-,j)$ is a
    conjugate pair of tunneling maps, and
    $\phi_1:\Sigma_1\rightarrow\Sigma_1$ is a diffeomorphism, then
    $(\phi_1^\ast v_+,\phi_1^\ast v_-,\phi_1^\ast j)$ is also a pair
    of conjugate tunneling maps with domain $(\phi_1^\ast\dot
    S,\phi_1^\ast j_1)$. This concludes the proof of Lemma
    \ref{lem:gauge_action}.}
\end{remark}

\vskip1cm
\section{Properties of the Folded Diagonal}
\label{sec:prop_of_Delta^Z}
In this section we study the properties of the folded diagonal by
looking at its deformation space. In short, we find that the
deformations of the folded diagonal at a pair of conjugate tunneling
maps restricted to $\sigma$ $(\hat v_+,\hat v_-)$ is given by the
graph of a function from sections of $\hat v_+^\ast TZ$ to sections of
$\hat v_-^\ast TZ$, which we will describe.

We only consider the case of $S^1$-invariant folds. They have the
advantage that the $\mathcal{H}$-holomorphic map equations are linear
and upper triangular, which simplifies the proofs. Moreover, in the
$S^1$-invariant case most of the proofs are constructive which
provides further intuitive insight into the properties of folded
holomorphic maps.

First we give the linearization of the equation for tunneling maps.
\begin{lemma}\label{lem:linearized_tunneling}
  At a solution $(v,j)$, the linearizations of equations
  (\ref{eq:H_F}) and (\ref{eq:H_L}) are
  \begin{eqnarray}
    D^{F}_{(v,j)}(\xi,h)
    &=&
    \pi_F\nabla^{0,1}\xi    +\frac12J\,\pi_F\,dv\,h=0
    \in\Omega^{0,1}(v^\ast F)
    \label{eq:lin_H_F}\\
    D^{L}_{(v,j)}(\xi,h) &=&d\left[
    d(\alpha(\xi))\circ j+v^\ast(\iota_{\xi} d\alpha)\circ j+
    v^\ast\alpha\circ h\right]=0 \in
    \Omega^2(v^\ast L)\label{eq:lin_H_L}
  \end{eqnarray}
  where and $\nabla$ is the Levi-Civita connection with respect to the
  metric $\tilde g|_Z=\alpha\otimes\alpha+g_F$ on $Z$ and
  \begin{eqnarray*}
    \nabla^{0,1}&=&\frac12\left\{\nabla+J\nabla\circ j\right\}.
  \end{eqnarray*}
\end{lemma}

\begin{proof}
  Let $v_t$ be a family of tunneling maps with complex structures
  $j_t$, $\xi=\frac{d}{dt}\big|_{t=0}v_t$ and
  $h=\frac{d}{dt}\big|_{t=0}j$. Recall that in the $S^1$-invariant
  case $\nabla R=0$, $\nabla J=0$ and $\nabla \pi_F=0$. Thus
  \begin{eqnarray*}
    D^{F}_{(v,j)}(\xi,h)
    &=&\left.\frac{d}{dt}\right|_{t=0}
    \frac{1}{2}\pi_F\left\{dv_t+J\,dv_t\,j_t\right\}\\
    &=&\frac{1}{2}\pi_F\left\{\nabla_t dv_t
      +J\nabla_t(dv_t)\circ j
      +J\,dv\,h\right\}_{t=0}\\
    &=&\pi_F\,\nabla^{0,1}\xi
      +\frac12J\,\pi_F\,dv\,h.
  \end{eqnarray*}
  
  For the second equation we compute
  \begin{eqnarray*}
    \left.\frac{d}{dt}\right|_{t=0} v_t^\ast\alpha\circ j_t
    =(v^\ast\mathcal{L}_\xi\alpha)\circ j+v^\ast\alpha\circ h
    =d(\alpha(\xi))\circ j+v^\ast(\iota_\xi d\alpha)\circ j+
    v^\ast\alpha\circ h.
  \end{eqnarray*}
\end{proof}
So the system of PDEs (\ref{eq:lin_H_F}) and (\ref{eq:lin_H_L}) is
upper triangular with respect to the splitting $TZ=L\oplus F$. This
greatly simplifies computations and the construction of examples.

Definition \ref{def:conjugate_tunneling} allows us to define several
maps.
\begin{definition}\label{def:F_isomorphism}
  Let $v_+$ and $v_-$ be conjugate tunneling maps. Then define the 
  $(J_+,J_-)$-linear vector bundle isomorphisms
  \begin{eqnarray*}
    A^F:v_+^\ast F\rightarrow v_-^\ast F,&\qquad 
    A^F=(\pi_F\,dv_-)\circ(\pi_F\,dv_+)^{-1}&\\
    A^E:v_+^\ast E\rightarrow v_-^\ast E\,&\qquad 
    A^E(\zeta_1\del_r+\zeta_2R)=\zeta_1\del_r-\zeta_2R.&
  \end{eqnarray*}
  $A^F$ is well defined at the zeros of $\pi_F\,dv_+$ by Equation
  (\ref{eq:conjugate_F}).
 
  By the vanishing condition on $\lambda$ along $Q$ we may define
  \begin{eqnarray}
    \label{eq:def_f}
    f:v_+^\ast F\rightarrow \setR,\qquad 
    f=\lambda\circ(\pi_F\,dv_+^{-1}).
  \end{eqnarray}
\end{definition}

\begin{remark}\label{rem:transversality_assumptions}
  In the following we will assume that, unless $\alpha$ is closed, at
  each point $q\in Q$ the map
  \begin{eqnarray}
    \label{eq:nabla_f}
    (\nabla f)(q):F_{v_\pm(q)}\rightarrow T^\ast_qS
  \end{eqnarray}
  is surjective and that the order of vanishing of $\pi_F\,dv_\pm$ is
  no greater than 1. We will deal with cases where these assumptions
  are not satisfied in the context of compactness of the moduli space.
\end{remark}

\begin{lemma}
  Let $(\xi_+,h)$ and $(\xi_-,h)$ be deformations of the conjugate
  tunneling maps $v_+$ and $v_-$
  The linearizations of the equations for conjugate tunneling
  maps (\ref{eq:conjugate_F}) through (\ref{eq:conjugate_marker}) are
  \begin{eqnarray}
    \label{eq:lin_conjugate_F}
    \pi_F\xi_+&=&\pi_F\xi_-\\
    d\left[\alpha(\xi_+)+\alpha(\xi_-)\right]\circ j
    &=&\lambda\circ jh
    -\left[v_+^\ast(\iota_{\xi_+}d\alpha)+v_-^\ast(\iota_{\xi_-}d\alpha)
      \right]\circ j
    \ \mathrm{on}\ T\sigma
    \label{eq:lin_conjugate_L}\\
    \alpha(\xi_+)(p)+\alpha(\xi_-)(p)&=&0.\label{eq:lin_conjugate_marker}
  \end{eqnarray}
\end{lemma}

\begin{proof}
  This follows immediately from the definitions and Remark
  \ref{rem:conjugate}.
\end{proof}

We need the following technical lemmas.
\begin{lemma}\label{lem:u-v-lambda}
  Let $(\Sigma,\sigma,j,u_+,u_-)$ be a folded holomorphic map with
  conjugate tunneling maps $(\dot S,j,v_+,v_-)$. Then for $\eta\in
  T_\sigma S=T_\sigma\Sigma$
  \begin{eqnarray}\label{eq:u-v-lambda}
    \lambda(\eta)=v_+^\ast\alpha(j\eta)+v_-^\ast\alpha(j\eta)
    =u_+^\ast\alpha(j\eta)+u_-^\ast\alpha(j\eta)
  \end{eqnarray}
\end{lemma}
\begin{proof}
  The first equality is just the definition of $\lambda$ (Equation
  \ref{eq:conjugate_L}). Recall that $\hat u_\pm=\hat v_\pm$ and both
  $du_\pm$ and $\pi_F\,dv_\pm$ are $(j,J_\pm)$-linear. Therefore
  \begin{eqnarray*}
    \pi_F\,du_\pm=\pi_F\,dv_\pm\qquad\mathrm{over}\ \sigma,
  \end{eqnarray*}
  and for $\eta\in T\sigma$ we have
  \begin{eqnarray*}
    u_\pm^\ast \alpha(\eta)=v_\pm^\ast\alpha(\eta),\qquad
    u_\pm^\ast\alpha(j\eta)=0,\qquad
    v_+^\ast\alpha(j\eta)=-v_-^\ast\alpha(j\eta)
  \end{eqnarray*}
  where the second equation follows from the fact that $\pi_E\,du_\pm$
  is $J_\pm$-linear and $\pi_K du_\pm|_{T\sigma}=0$ and the last
  equation is just Equation (\ref{eq:conjugate_L}). Equation
  (\ref{eq:u-v-lambda}) follows.
\end{proof}

\begin{lemma}\label{lem:eta_formula}
  Let $v$ be a tunneling map, $\eta$ a section of $TS$ vanishing at
  the puncture, and $h$ a deformation of complex structure $j$ such
  that
  \begin{eqnarray*}
    \nabla^{0,1}\eta+\frac12jh=0.
  \end{eqnarray*}
  Then
  \begin{eqnarray}\label{eq:eta_formula}
    d(v^\ast\alpha(j\eta))
    -d(v^\ast\alpha(\eta))\circ j
    =(\iota_\eta v^\ast d\alpha)\circ j +v^\ast\alpha\circ h
  \end{eqnarray}
  and consequently, if $v_+$ and $v_-$ are conjugate tunneling maps,
  \begin{eqnarray}\label{eq:lambda_formula}
    d[\lambda(\eta)]+d[\lambda(j\eta)]\circ j=
    (\iota_\eta [v_+^\ast d\alpha+v_-^\ast d\alpha])\circ j 
    -\lambda\circ jh.
  \end{eqnarray}
\end{lemma}

\begin{proof}
  Near some point $z\in \dot S$, let $X$ be a section of $TS$ with
  $\nabla X(z)=0$.
  \begin{eqnarray*}
    d(v^\ast\alpha(\eta))\circ j(X)
    &=&\mathcal{L}_{jX}(v^\ast\alpha(\eta))\\
    &=&\mathcal{L}_{jX}(v^\ast\alpha)(\eta)
    +v^\ast\alpha(\nabla_{jX}\eta)\\
    &=&\mathcal{L}_\eta(v^\ast\alpha(jX))+v^\ast d\alpha(jX,\eta)
    +v^\ast\alpha(j\nabla_X\eta- h(X))\\
    &=&\mathcal{L}_\eta(v^\ast\alpha\circ j)(X)-v^\ast d\alpha(\eta,jX)
    -v^\ast\alpha\circ h(X)\\
    &=&[d(v^\ast\alpha(j\eta))-(\iota_\eta v^\ast
    d\alpha)\circ j -v^\ast\alpha\circ h](X).
  \end{eqnarray*}
\end{proof}

\begin{lemma}\label{lem:f_formula}
  Let $(v_+,v_-,j)$ be a pair of conjugate tunneling maps and, $\xi$ a
  section of $v_+^\ast TZ$ and $h$ a
  deformation of complex structure $j$ such that
  \begin{eqnarray*}
    D^F_{(v_+,j)}(\xi,h)=0.
  \end{eqnarray*}
  Then, with $\chi=\pi_F\xi$,
  \begin{eqnarray}\label{eq:f_formula}
    d[f(\chi)]+d[f(J\chi)]\circ j
    =[v_+^\ast(\iota_\chi d\alpha)
    +v_-^\ast(\iota_{A^F(\chi)} d\alpha)]\circ j
    -\lambda\circ jh.
  \end{eqnarray}
\end{lemma}

\begin{proof}
  Note that $f$ is a smooth function, so it suffices to show this on
  the dense subset $\ddot S=\dot S\setminus Q\subset\dot S$.  On
  $\ddot S$ we may write $\chi=\pi_F\,dv_+(\eta)$.

  Note that on $\ddot S$
  \begin{eqnarray*}
    D^F_{(v_+,j)}(\xi,h)=0\qquad
    \mbox{if and only if}\quad
    \nabla^{0,1}\eta+\frac12jh=0.
  \end{eqnarray*}
  Using Lemma \ref{lem:eta_formula} and the fact that
  $\iota_Rd\alpha=0$
  \begin{eqnarray*}
    d[f(\chi)]+d[f(J\chi)]\circ j
    &=&d[\lambda(\eta)]+d[\lambda(j\eta)]\circ j\\
    &=&(\iota_\eta[v_+^\ast d\alpha+v_-^\ast d\alpha])\circ j
    -\lambda\circ jh\\
    &=&[v_+^\ast(\iota_\chi d\alpha)
    +v_-^\ast (\iota_{A^F(\chi)}d\alpha)]\circ j -\lambda\circ jh
  \end{eqnarray*}
  so Equation \ref{eq:f_formula} holds on $\ddot S$ and by smoothness
  of $f$ it holds on $\dot S$.
\end{proof}

\subsection{Deformations of the Folded Diagonal}
Using the results of the previous section we describe the deformations
of the folded diagonal. We show that they are given by the graph of a
pseudo-differential operator $B^Z:\Gamma(\hat u_+^\ast TZ)\rightarrow
\Gamma(\hat u_-^\ast TZ)$. 

\begin{theorem}\label{thm:deformations}
  Given conjugate tunneling maps $(v_+,v_-,j)$, a section $\hat\xi$ of
  $\hat v_+^\ast TZ$.
  Then there exists an extension $\xi_+$ of $\hat \xi$ to $v_+^\ast
  TZ$ a section $\xi_-$ of $v_-^\ast TZ$ and a deformation $h$ of the
  complex structure $j$ so that $(\xi_\pm,h)$ satisfy Equations
  (\ref{eq:lin_H_F}) and (\ref{eq:lin_H_L}) and $(\xi_+,\xi_-,h)$
  satisfies Equations (\ref{eq:lin_conjugate_F}) through
  (\ref{eq:lin_conjugate_marker}) and preserving the condition on
  $\lambda$ at the non-immersion points.

  Moreover the restriction $\xi_-|_\sigma$ is uniquely determined by
  $\hat\xi$.
\end{theorem}

\begin{proof}
  We split this argument into the two cases when $\alpha$ is closed
  and when $\alpha$ is not closed.

  If $\alpha$ is not closed, recall Definition \ref{eq:def_f} of the
  function $f:v_\pm^\ast F\rightarrow \setR$. Let
  $g:S\rightarrow\setR$ be the unique harmonic function vanishing at
  the puncture and satisfying von Neumann boundary conditions
  \begin{eqnarray}\label{eq:def_g}
    dg\circ j=d[f(\pi_F\hat\xi)]\qquad\mathrm{on}\ T\sigma.
  \end{eqnarray}
  
  Let $\xi$ be any extension of $\hat\xi$ to $\dot S$ that vanishes
  in a neighborhood of the puncture, satisfies 
  \begin{eqnarray}\label{eq:nabla_f_Q}
    (\nabla f)(\pi_F\xi)=dg\circ j\qquad\mathrm{on}\ Q     
  \end{eqnarray}
  and $\pi_F\nabla^{0,1}\xi=0$ in a
  neighborhood of $Q$. This is well defined by the assumption that
  $f|_Q:F_{v(Q)}\rightarrow T^\ast_QS$ is surjective (see Remark
  \ref{rem:transversality_assumptions}).

  Then we may set
  \begin{eqnarray*}
    h=2 (\pi_F\,dv_+)^{-1}(J\pi_F\nabla^{0,1}\xi)
  \end{eqnarray*}
  which is $j$-antilinear, and consequently $D^F_{(v_+,j)}(\xi,h)=0$.
  Now let $\zeta:S\rightarrow\setR$ be the unique solution to the
  Dirichlet problem
  \begin{eqnarray*}
    d(d\zeta\circ j)&=&-d\left[
    d(\alpha(\xi))\circ j+v_+^\ast(\iota_{\xi} d\alpha)\circ j+
    v_+^\ast\alpha\circ h\right],\\
    \zeta|_\sigma&=&0.
  \end{eqnarray*}
  Then $\xi_+=\xi+\zeta\cdot R$ satisfies 
  \begin{eqnarray*}
    D^F_{(v_+,j)}(\xi_+,h)&=&D^F_{(v_+,j)}(\xi,h)=0,\\
    D^L_{(v_+,h)}(\xi_+,h)&=&D^L_{(v_+,j)}(\xi,h)+d(d\zeta\circ j)=0,\\
    \xi_+|_\sigma&=&\hat\xi.
  \end{eqnarray*}
  
  Now define the section $\xi_-$ of $v_-^\ast TZ$ by
  \begin{eqnarray}\label{eq:xi_-}
    \xi_-=A^F(\pi_F\xi_+)-(f(J\pi_F\xi_+)+\alpha(\xi_+)+g) R.
  \end{eqnarray}
  Then
  \begin{eqnarray*}
    D^F_{(v_-,j)}(\xi_-,h)=D^F_{(v_+,j)}(\xi_+,h)=0
  \end{eqnarray*}
  and
  \begin{eqnarray*}
    D^L_{(v_-,j)}(\xi_-,h)
    &=&d\left[
    d[\alpha(\xi_-)]\circ j+v_-^\ast(\iota_{\xi_-} d\alpha)\circ j+
    v_-^\ast\alpha\circ h\right]\\
    &=&d\left[
    -d[f(J\pi_F\xi_+)+\alpha(\xi_+)+g]\circ j+v_-^\ast(\iota_{\xi_-}
    d\alpha)\circ j -\lambda\circ jh
    -v_+^\ast\alpha\circ h\right]\\
    &=&d\left[-d[f(J\pi_F\xi_+)]\circ j
      +v_+^\ast(\iota_{\xi_+} d\alpha)\circ j
      +v_-^\ast(\iota_{\xi_-} d\alpha)\circ j 
      -\lambda\circ jh\right] \\
    &=&d\left[d(f(\pi_F\xi_+))\right]\\
    &=&0
  \end{eqnarray*}
  where we used Equation \ref{eq:f_formula} in the second last
  step. 
  
  Now we need to show that $(\xi_+,\xi_-,h)$ is a deformation of
  conjugate tunneling maps. By construction
  $\pi_F\xi_-=A^F(\pi_F\xi_+)$, verifying Equation
  (\ref{eq:lin_conjugate_F}). For Equation (\ref{eq:lin_conjugate_L})
  we compute on $T\sigma$
  \begin{eqnarray*}
    d[\alpha(\xi_+)+\alpha(\xi_-)]\circ j
    &=&-d[f(J\pi_F\xi_+)+g]\circ j\\
    &=&-d[f(J\pi_F\xi_+)]\circ j-dg\circ j\\
    &=&d[f(\pi_F\xi_+)]-dg\circ j
    -[v_+^\ast(\iota_{\xi_+} d\alpha)
    +v_-^\ast(\iota_{\xi_-} d\alpha)]\circ j+\lambda\circ jh\\
    &=&d[f(\pi_F\hat\xi)]-dg\circ j
    -[v_+^\ast(\iota_{\xi_+} d\alpha)
    +v_-^\ast(\iota_{\xi_-} d\alpha)]\circ j+\lambda\circ jh\\
    &=&\lambda\circ jh
    -[v_+^\ast(\iota_{\xi_+} d\alpha)
    +v_-^\ast(\iota_{\xi_-} d\alpha)]\circ j
  \end{eqnarray*}
  where we used Equations (\ref{eq:f_formula}) and (\ref{eq:def_g}).
  For Equation (\ref{eq:lin_conjugate_marker}) we compute
  \begin{eqnarray*}
    \alpha(\xi_+)(p)+\alpha(\xi_-)(p)=f(J\pi_F\xi_+(p))+g(p)
    =f(J\pi_F\xi(p))=0.
  \end{eqnarray*}
  
  We are left to verify that the condition on the order of vanishing
  of $\lambda$ along $Q$ is preserved. Linearizing
  $\pi_F\,dv_\pm|_Q=0$ gives
  \begin{eqnarray*}
    \pi_F\nabla_\nu \pi_F\,dv_\pm+\pi_F\nabla(\pi_F\xi_\pm)=0
    \qquad\mathrm{at\ all}\ q\in Q
  \end{eqnarray*}
  where $\nu$ is the deformation of $q$. Note that
  $\pi_F\nabla_\nu(\pi_F\,dv_\pm)=\pi_F\nabla(\pi_F\,dv_\pm(\nu))$ at
  $q$, so we conclude that
  \begin{eqnarray}\label{eq:nu}
   \pi_F\nabla\left(\pi_F(dv_\pm(\nu)+\xi_\pm)\right)=0
    \qquad\mathrm{at\ all}\ q\in Q.
  \end{eqnarray}

  Linearizing the condition that $\lambda$ vanishes on $Q$ gives
  \begin{eqnarray}\label{eq:lin_lambda_Q}
    \nabla_\nu\lambda
    +[v_+^\ast(\iota_{\xi_+}d\alpha)
    +v_-^\ast(\iota_{\xi_-}d\alpha)]\circ j 
    +d[\alpha(\xi_+)+\alpha(\xi_-)]\circ j -\lambda\circ jh=0.
  \end{eqnarray}
  Using Equation (\ref{eq:f_formula}) we see that
  \begin{eqnarray*}
    &&[v_+^\ast(\iota_{\xi_+}d\alpha)
    +v_-^\ast(\iota_{\xi_-}d\alpha)]\circ j
    +d[\alpha(\xi_+)+\alpha(\xi_-)]\circ j
    -\lambda\circ jh\\
    &=&[v_+^\ast(\iota_{\xi_+}d\alpha)
    +v_-^\ast(\iota_{\xi_-}d\alpha)]\circ j
    -d[f(J\pi_F\xi_\pm)]\circ j-dg\circ j
    -\lambda\circ jh\\
    &=&d[f(\pi_F\xi_\pm)]-dg\circ j.
  \end{eqnarray*}
  Using that $\nabla_\nu\lambda=d(\lambda(\nu))=d[f(dv_\pm(\nu))]$ on
  $Q$ we see that on $Q$
  \begin{eqnarray*}
    \nabla_\nu\lambda+d[f(\pi_F\xi_\pm)]-dg\circ j
    &=&\nabla\left(f[\pi_F(\xi_\pm+dv_\pm(\nu))]\right)-dg\circ j\\
    &=&(\nabla f)\left(\pi_F[\xi_\pm+dv_\pm(\nu)]\right)-dg\circ j\\
    &=&(\nabla f)(\pi_F\xi_\pm)-dg\circ j\\
    &=&0,
  \end{eqnarray*}
  where we used Equation (\ref{eq:nu}) in the third last step and
  Equation (\ref{eq:nabla_f_Q}) in the last step. Combining this with
  the previous computation verifies equation (\ref{eq:lin_lambda_Q}).
  
  This establishes existence. To see uniqueness, suppose that
  $(\xi_+',h')$ also satisfies $D_{(v_+,j)}(\xi_+',h')=0$ and
  $\xi_+'|_\sigma=\hat\xi$. Then $(\xi''=\xi_+-\xi_+',h''=h-h')$
  satisfies
  \begin{eqnarray*}
    D(\xi'',h'')=0,\qquad \xi''|_\sigma=0.
  \end{eqnarray*}
  The above argument shows that given $(\xi_+,h)$ there is a unique
  pair $(\xi_-,h)$ satisfying Equations (\ref{eq:lin_H_F}) and
  (\ref{eq:lin_H_L}) so that $(\xi_+,\xi_-,h,\mu)$ satisfy Equations
  (\ref{eq:lin_conjugate_F}) through (\ref{eq:lin_conjugate_marker}).
  In order that the condition on the vanishing of $\lambda$ is
  preserved we need that, with the above notation,
  $f(J\pi_F\xi_+)=dg\circ j$ on $Q$. Thus we conclude that
  \begin{eqnarray*}
    \pi_F\xi''=0\qquad\mathrm{on}\ Q,
  \end{eqnarray*}
  so there exists a vector field $\eta$ on $S$ vanishing on $\sigma$
  so that $\pi_F\xi''=\pi_F\,dv_+(\eta)$. Since $\pi_F\,dv_+$ is
  $(j,J)$ linear and $\nabla J=0$ we conclude that
  $\nabla^{0,1}\eta+\frac12jh''=0$. The function
  $\zeta''=\alpha(\xi''-dv_+(\eta))$ vanishes on $\sigma$ and
  satisfies
  \begin{eqnarray*}
    d(d\zeta''\circ j)=0
  \end{eqnarray*}
  so $\zeta''=0$ and $\xi''=dv_+(\eta)$. Thus $(\xi_+,h)$ and
  $(\xi_+',h')$ differ by an infinitesimal gauge transformation $\eta$.
  By Remark \ref{rem:gauge_invariant} this gauge ambiguity preserves
  the conjugacy condition and does not change the boundary values
  $\xi_\pm|_\sigma$.
  
  Now consider the case when $d\alpha=0$. In this case $Z$ is a
  trivial $S^1$-bundle over $V$ and the two Equations
  (\ref{eq:lin_H_F}) and (\ref{eq:lin_H_L}) decouple and we may solve
  each one separately. This implies that
  $v_+^\ast\alpha+v_-^\ast\alpha\equiv 0$, so $\lambda\equiv 0$. As
  before, for given boundary values $\hat\xi\in\Gamma(\hat v_+^\ast
  TZ)$ there exist extensions $\xi_+$ to sections of $v_+^\ast TZ$ and
  deformation $h$ of complex structure $j$ so that Equations
  (\ref{eq:lin_H_F}) and (\ref{eq:lin_H_L}) are satisfied. Then
  \begin{eqnarray*}
    \xi_-=A^F(\pi_F\xi_+)-\alpha(\xi_+)R
  \end{eqnarray*}
  is the unique section of $v_-^\ast TZ$ satisfying Equations
  (\ref{eq:lin_H_F}), (\ref{eq:lin_H_L}) and
  (\ref{eq:lin_conjugate_F}) through (\ref{eq:lin_conjugate_marker}).
  But $\xi_-|_\sigma$ only depends on $\xi_+|_\sigma$, so all possible
  choices of $\xi_+$ lead to the same boundary values
  $\hat\xi_-=\xi_-|_\sigma$, concluding the proof of the theorem.
\end{proof}

\begin{definition}
  Given conjugate tunneling maps $(v_+,v_-)$, let $d\Delta^Z$ be the
  space of deformations of conjugate tunneling maps, restricted to
  $\sigma$, i.e.
  \begin{eqnarray*}
    d\Delta^Z&=&\{(\hat\xi_+,\hat\xi_-)\in \Gamma(\hat v_+^\ast
    TZ\oplus \hat v_-^\ast TZ)\big|\exists\, 
    \xi_\pm\in\Gamma(v_\pm^\ast TZ), h\in
    T_j\mathcal{J}(S)\ \mathrm{with}\ \xi_\pm|_\sigma=\hat\xi_\pm\\
    &&\mbox{ satisfying Equations }D_{(v_\pm,j)}(\xi_\pm,h)=0,\mbox{\ 
      (\ref{eq:lin_conjugate_F}), 
      (\ref{eq:lin_conjugate_L}),
      (\ref{eq:lin_conjugate_marker}),}\\
    &&\mbox{ and the vanishing condition on }\lambda\mbox{ along }Q.
      \}
  \end{eqnarray*}
\end{definition}

By the Definition \ref{def:folded_maps} of folded maps we identify
$\sigma=\sigma_0=\sigma_1$ and therefore we may identify $\hat
v_\pm^\ast TZ$ with $\hat u_\pm^\ast TZ$, when $(u_+,u_-)$ is a folded
map with tunneling maps $(v_+,v_-)$. Thus we may also view space
of deformations of folded diagonal $d\Delta^Z$ as a subset of
\begin{eqnarray*}
  \Gamma(\hat u_+^\ast TZ\oplus\hat u_-^\ast TZ)\subset \Gamma(\hat
  u_+^\ast TX\oplus \hat u_-^\ast TX).
\end{eqnarray*}
Note that if $(\hat\xi_+,\hat\xi_-)\in d\Delta^Z$, and $\eta\in
T_\id\Diff(\Sigma_1,\sigma)$ is an infinitesimal gauge
transformation that is tangent to the domain fold $\sigma$, then also
$(\hat\xi_++du_+(\eta),\hat\xi_-+du_-(\eta))\in d\Delta^Z$. This
defines an action
\begin{eqnarray}
  T_\id\Diff(\Sigma,\sigma)\times d\Delta^Z\rightarrow d\Delta^Z
  \label{eq:gauge_action}.
\end{eqnarray}

It will prove convenient to extend the definition of deformations of
folded diagonal to sections of $\hat u_+^\ast TX\oplus \hat u_-^\ast
TX$ in such a way that the action (\ref{eq:gauge_action}) extends to
all infinitesimal gauge transformations $T_\id\Diff(\Sigma)$ of the
map domain, including the ones that move the domain fold. This will
greatly simplify taking the quotient by the gauge action later. The
following definitions and lemmas facilitate this.

\begin{definition}
  Let $(u_+,u_-)$ be a folded holomorphic map with conjugate tunneling
  maps $(v_+,v_-)$.
  
  Let $H_\pm=u_\pm^\ast TX$ and set $\hat H_\pm=H_\pm|_\sigma$.  We
  define the subbundles of $\hat H_\pm$
  \begin{eqnarray*}
    F_\pm=\hat u_\pm^\ast F,\qquad
    E_\pm=\hat u_\pm^\ast E.
  \end{eqnarray*}
  
  Recall the $(J_+,J_-)$-linear operators $A^F:F_+\rightarrow F_-$ and
  $A^E:E_+\rightarrow E_-$ defined in Definition \ref{def:F_isomorphism}.
  Let
  \begin{eqnarray*}
    \tilde Q&:&\Gamma(E_-)\rightarrow \Gamma(v_-^\ast E)\\
  \end{eqnarray*}
  be the operator so that for $\zeta^K,\zeta^L:\sigma\rightarrow
  \setR$ we have that $f\del_r+g\,R=\tilde Q(\zeta^K\del_r+\zeta^L R)$
  satisfies
  \begin{eqnarray*}
    df\circ j+dg&\in&\mathcal{H}_n(S)\\
    f|_\sigma&=&\zeta^K\\
    g(p)&=&0
  \end{eqnarray*}
  and define
  \begin{eqnarray*}
    D&:&\Gamma(E_-)\rightarrow \Gamma(E_-),\qquad 
    D(\zeta^K\cdot\del_r+\zeta^L\cdot R)
    =\zeta^K\cdot\del_r-\zeta^L\cdot R\\
    Q&:&\Gamma(E_-)\rightarrow \Gamma(E_-),\qquad
    Q(\zeta)=\tilde Q(\zeta)|_\sigma\\
    C&:&\Gamma(E_+)\rightarrow\Gamma(E_-),\qquad
    C=-D\circ Q\circ A^E\\
    f_\setC&:&F_+\rightarrow E_+,\qquad
    f_\setC(\chi)=f(\chi)\otimes\del_r-f(J\chi)\otimes R.
  \end{eqnarray*}
  Note that $D$ is $J_-$ anti-linear and $f_\setC$ is $J_+$ linear.
\end{definition}

\begin{lemma}
  $C$ is a pseudo-differential operator of order zero with principal
  symbol
  \begin{eqnarray}
    \label{eq:c}
    c(s,\zeta)=-(\id-iJ_-)\pi_K\, A^E\zeta.
  \end{eqnarray}
\end{lemma}
\begin{proof}
  We claim that $Q$ is a pseudo-differential operator of order zero
  with principal symbol
  \begin{eqnarray*}
    q(s,\zeta)=\frac{1}{2}(\id +iJ_-)(\zeta+D\zeta)=(\id +iJ_-)\pi_K\zeta.
  \end{eqnarray*}
  To see this we note that $Q$ can be written as
  \begin{eqnarray*}
    Q=P+D\circ (1-P)
  \end{eqnarray*}
  where $P$ is the Calder\'on projector onto the space of Cauchy data
  of the elliptic operator $\tilde Q$. Then recall that $\sigma$ has
  the opposite orientation from $\del S$, so the principal symbol $p$
  of $P$ is just the projection $p=\frac12(\id +iJ_-)$ onto the $J_-$
  anti-holomorphic subspace of $E_-$. For a detailed discussion of
  this see \cite{seeley} or \cite{booss_wojchiechowski} in the case of
  Cauchy-Riemann operators and \cite{liviu}, \S 6.2 for the family
  version, interpreting $Q$ as a family of Cauchy-Riemann operators
  with parameter space $\mathcal{H}_n(S)$.
  
  Thus
  \begin{eqnarray*}
    c=-D\circ q\circ A^E
    =-D(\id+iJ_-)\pi_K\,A^E
    =-(\id-iJ_-)D\pi_K\, A^E
    =-(\id-iJ_-)\pi_K\, A^E.
  \end{eqnarray*}
\end{proof}

\begin{definition}\label{def:B}
  Let $(\Sigma,\sigma,j,u_+,u_-)$ be a folded holomorphic map.
  Using the splitting $\hat H_\pm =F_\pm\oplus E_\pm$
  we define the linear map on sections
  \begin{eqnarray}
    B&:&\Gamma(\hat H_+)\rightarrow \Gamma(\hat H_-)\nonumber\\
    B(\xi^F\oplus \xi^E)&=&A^F(\xi^F)\oplus A^E(\xi^E-f_\setC(\xi^F))
    +C((1-a)\xi^E-f_\setC(\xi^F)),
    \label{eq:B} 
  \end{eqnarray}
  where $a$ is the gap-function from Equation (\ref{eq:a_tau}).
\end{definition}

\begin{lemma}\label{lem:B gives dDeltaZ}
  The deformations of the folded diagonal are given by the graph of
  the operator $B$ defined in (\ref{eq:B}), restricted to
  $\hat u_+^\ast TZ$, i.e.
  \begin{eqnarray*}
    d\Delta^Z=\graph(B|_{\hat u_+^\ast TZ}).
  \end{eqnarray*}
\end{lemma}

\begin{proof}
  For $\hat\xi\in\Gamma(\hat u_+^\ast TZ)$
  \begin{eqnarray*}
    B(\hat \xi)
    &=&A^F[\pi_F\hat\xi]+A^E[\pi_E\hat\xi-f_\setC(\pi_F\hat\xi)]
    +C[(1-a)\pi_E\hat\xi-f_\setC(\pi_F\hat\xi)]\\
    &=&A^F[\pi_F\hat\xi]-[\alpha(\hat\xi)+f(J\pi_F\hat\xi)] R
    -f(\pi_F\hat\xi)\del_r+D\circ Q[f(\pi_F\hat\xi)\del_r]\\
    &=&A^F[\pi_F\hat\xi]-[\alpha(\hat\xi)+f(J\pi_F\hat\xi)] R
    -\pi_L Q[f(\pi_F\hat\xi)\del_r]
  \end{eqnarray*}
  Then $Q(f(\pi_F\xi)\del_r)$ is exactly the function $g$ from
  Equation (\ref{eq:def_g}) from the proof of Theorem
  \ref{thm:deformations}, so $B(\hat \xi)$ equals $\xi_-|_\sigma$ from
  Equation (\ref{eq:xi_-}).
\end{proof}

Using $B$ we extend the definition of the folded diagonal to a
subset of $\Gamma(\hat u_+^\ast TX\oplus\hat u_-^\ast TX)$ in the
obvious way:
\begin{definition}\label{def:extended deformations}
  The space of extended deformations $d\Delta^X$ of the folded
  diagonal is
  \begin{eqnarray*}
    d\Delta^X&=&\graph(B).
  \end{eqnarray*}
\end{definition}

The following result motivates the definition of $B$ and $d\Delta^X$.
\begin{lemma}\label{lem:B_gauge_invariant}
  The space of extended deformations $d\Delta^X$ of the folded
  diagonal is invariant under the full infinitesimal gauge group of
  $\Sigma$ (not just the subgroup that preserves $\sigma$ as a set).
\end{lemma}
\begin{proof}
  Let $\eta$ be a section of $T_\sigma S=T_\sigma\Sigma$. Then
  \begin{eqnarray*}
    f(\pi_F\,du_+(\eta))
    &=&f(\pi_F\,dv_+(\eta))
    =\lambda(\eta)\qquad\mathrm{and}\\
    (1-a)\pi_K\,du_+(\eta)
    &=&(1-a)u_+^\ast\alpha(j\eta)\del_r
    =u_+^\ast\alpha(j\eta)+u_-^\ast\alpha(j\eta)
    =\lambda(\eta)\del_r
  \end{eqnarray*}
  and $C[(1-a)\pi_E\,du_+(\eta)-f_\setC(\pi_F\,du_+(\eta))]=C(0)=0$.
  Then
  \begin{eqnarray*}
    B(du_+(\eta))
    &=&A^F(\pi_F\,du_+(\eta))
    +A^E[\pi_E\,du_+(\eta)-f_\setC(\pi_F\,du_+(\eta))]\\
    &=&A^F(\pi_F\,dv_+(\eta))
    +A^E[(u_+^\ast\alpha(j\eta)-\lambda(\eta))\del_r
    +(u_+^\ast\alpha(\eta)+\lambda(j\eta)) R]\\
    &=&\pi_F\,dv_-(\eta)
    +(u_+^\ast\alpha(j\eta)-\lambda(\eta))\del_r
    -(u_+^\ast\alpha(\eta)+\lambda(j\eta)) R\\
    &=&\pi_F\,du_-(\eta)
    -u_-^\ast\alpha(j\eta)\del_r
    +u_-^\ast\alpha(\eta) R\\
    &=&\pi_F\,du_-(\eta)
    +\pi_E\,du_-(\eta)\\
    &=&du_-(\eta).
  \end{eqnarray*}
\end{proof}

\vskip1cm
\section{The Moduli Space of Folded Holomorphic Maps}
Now we come to the result that justifies the definitions and lemmas
pertaining to tunneling maps. We show that they give elliptic boundary
values.

Let $(u_+,u_-,j)$ be a folded holomorphic map with domain fold
$\sigma\subset \Sigma$.
\begin{theorem}\label{thm:elliptic boundary conditions}
  Assume that the map $(u_+,u_-)$ is transverse to the fold, so
  $\sigma$ is a manifold. Then the map
  \begin{eqnarray*}
    R:\Gamma(\hat H_+)\oplus \Gamma(\hat H_-)\rightarrow 
    \Gamma(\hat H_-),\qquad
    \xi_+\oplus\xi_-\mapsto \xi_--B(\xi_+)
  \end{eqnarray*}
  poses elliptic boundary conditions for the folded holomorphic map
  $(u_+,u_-,j)$, i.e. the principal symbol $r$ of $R$ restricted to
  the range of the principal symbol $p$ of the Calder\'on projector
  $P$ for the complexified Cauchy-Riemann operator $D_{u_+}\times
  D_{u_-}$ on $H_+\otimes\setC\times H_-\otimes\setC$
  \begin{eqnarray*}
    r|_{\range(p)}\rightarrow \hat H_-\otimes\setC
  \end{eqnarray*}
  is an isomorphism.
\end{theorem}

\begin{proof}
  Note that the principal symbol $r$ is given by $r(u,v)=v-b(u)$,
  where $b$ is the principal symbol of $B$.
  
  Let $\hat H_\pm^\setC=\hat H_\pm\otimes \setC$ denote the
  complexification of $\hat H_\pm$ and let $\hat H_\pm'$ ($\hat
  H_\pm''$) denote the $(i,J_\pm)$-linear (antilinear) subspace of
  $\hat H_\pm^\setC$. Recall that $\sigma=\del\Sigma_+=-\del\Sigma_-$
  inherits the orientation from $\Sigma_+$. Then
  \begin{eqnarray*}
    p&:&\hat H_+^\setC\oplus \hat H_-^\setC\rightarrow \hat
    H_+'\oplus \hat H_-''\\
    p(v,w)&=&\frac12(\id-i J_+)\oplus\frac12(\id+iJ_-)
  \end{eqnarray*}
  is the projection onto the $(i,J_+\oplus -J_-)$ linear subspace.
  
  Note that $c=\frac12(\id-iJ_-)c$ equals its projection onto
  $E'$. Evaluating $c$ on an element $w\in E_+'$ gives
  \begin{eqnarray*}
    c(w)
    &=&-(\id-iJ_-)\pi_K\,A^E(w)
    =-A^E[(\id-iJ_+)\pi_K w]\\
    &=&-A^E[\pi_K w-\pi_L(iJ_+w)]
    =-A^E\left(\pi_K w+\pi_L w\right)\\
    &=&-A^E w.
  \end{eqnarray*}

  With the notation from Definition \ref{def:B} and $w=w^F \oplus
  w^E\in F_+'\oplus E_+'$
  \begin{eqnarray*}
    b(w)&=&A^F(w^F)+A^E(w^E-f_\setC(w^F))
    +c\left((1-a)w^E-f_\setC(w^F)\right)\\
    &=&A^F(w^F)+A^E\left(w^E-f_\setC(w^F)
      -(1-a)w^E+f_\setC(w^F)\right)\\
    &=&A^F(w^F)+a\,A^E(w^E).
  \end{eqnarray*}
  Now suppose that $(w,z)\in \range(p)=\hat H_+'\oplus \hat H_-''$ and
  $r(w,z)=0$. Then necessarily $w^F=0$ and $w^E=0$ as $a>0$,
  and therefore also $z=0$. We conclude that $r$ is an isomorphism.
\end{proof}

Now standard theory shows that the linearized operator is Fredholm.
For this next theorem we fix the complex structure on the domain
$\Sigma=\Sigma_0$.
\begin{theorem}\label{thm:fredholm}
For any $s>\frac12$, the operator
\begin{eqnarray*}
  D^s_B&:&H^s(\Sigma_+,H_+)\times H^s(\Sigma_-,H_-)\rightarrow\\
  &&H^{s-1}(\Sigma_+,\Lambda^{0,1}T^\ast\Sigma_+)\times
  H^{s-1}(\Sigma_-,\Lambda^{0,1}T^\ast\Sigma_-)\times
  H^{s-\frac12}(\del\Sigma_-,\hat H_-)\\
  &&(\xi_+,\xi_-)\mapsto
  \left(D_{u_+}\xi_+,D_{u_-}\xi_-,R(\hat\xi_+,\hat\xi_-)\right)
\end{eqnarray*}
is Fredholm with real Fredholm index
\begin{eqnarray*}
  index(D^s_B)=\mu(H_+,F_+)+\mu(H_-,F_-)+2\chi(\Sigma)
\end{eqnarray*}
where $F_\pm=K\oplus \pi_F\,du_\pm(T\sigma)$ are totally real
subbundles and $\mu$ is the Maslov index.

The kernel of $D^s_B$ is independent of choice of $s>\frac12$
and consists only of smooth solutions.
\end{theorem}

\begin{proof}
  The Fredholm properties and smoothness of solutions are a direct
  application of Theorem 19.1 and 20.8 of \cite{booss_wojchiechowski},
  using Theorem \ref{thm:elliptic boundary conditions} above. To see
  the index formula, define the homotopy of boundary conditions
  \begin{eqnarray*}
    B_t(\xi^F\oplus\xi^E)=A^F(\xi^F)+\,A^E(\xi^E-f_\setC(\xi^F))+t\cdot
    C\left((1-a)\xi^E -f_\setC(\xi^F)\right).
  \end{eqnarray*}
  Following the arguments as in the proof of Theorem \ref{thm:elliptic
    boundary conditions} we see that the symbol $b_t$ of $B_t$
  satisfies on $w=w^E+w^F\in E_+'\oplus F_+'$
  \begin{eqnarray*}
    b_t(w)
    &=&A^F(w^F)+A^E[(w^E-f_\setC(w^F)]-t\,A^E[(1-a)w^E-f_\setC(w^F)]\\
    &=&A^F(w^F)+A^E[(1+at-t)w^E-(1-t)f_\setC(w^F)].
  \end{eqnarray*}
  Thus if $b_t(w)=0$, then $w^F=0$ and $(1+at-t) A^E(w^E)=0$.  Using
  that $a>0$ we see that $1+at-t>0$ for $t\in[0,1]$ and we conclude
  that $w^E=0$ for all $t\in[0,1]$. Thus each member in this family
  gives elliptic boundary conditions.
  
  $B_0$ is the $(J_+,J_-)$-linear bundle isomorphism with
  \begin{eqnarray*}
    B_0(K)=K,\qquad
    B_0(\pi_F\,du_+(T\sigma))=\pi_F\,du_-(T\sigma),
  \end{eqnarray*}
  so $B_0(F_+)=F_-$.
  
  Then $\graph(B_0)$ has the same Maslov index as $(F_+,F_-)$, as can
  be seen by the homotopy of totally real subspaces
  \begin{eqnarray*}
    \Lambda_t=\{(u+(1-t)J_+v,B_0[(1-t-J_+t)((1-t)u+J_+v)])|u,v\in F_+\}
  \end{eqnarray*}
  satisfying
  \begin{eqnarray*}
    \Lambda_0=\graph(B_0),\qquad
    \Lambda_1=(F_+,B_0(F_+))=(F_+,F_-).
  \end{eqnarray*}
  
  In conclusion we have that the index for the boundary value problem
  given by $B=B_1$ is the same as the one for $B_0$, which in turn can
  be computed by $(F_+,F_-)$.
\end{proof}

To visualize the construction and the results up to here consider the
following. As seen in Section \ref{sec:motivating_example}, the
diagonal in $\Map(\sigma,Z)\times\Map(\sigma,Z)$ does not yield
elliptic boundary conditions.  In the language of \cite{liviu} we may
say that given holomorphic maps $(u_+,u_-)$ with
$u_+|_\sigma=u_-|_\sigma$, the subspace
\begin{eqnarray*}
  \{(\hat\xi_+,\hat\xi_-)\in L^2(\sigma,\hat u_+^\ast TX)\times
  L^2(\sigma,\hat u_-^\ast TX)|\exists\,
  \xi_\pm\in \ker(D_{u_\pm}),\ \xi_\pm|_\sigma=\hat\xi_\pm\}
\end{eqnarray*}
is not a Fredholm pair. But the folded diagonal gives elliptic
boundary conditions, or given a folded holomorphic map $(u_+,u_-)$ the
subspace
\begin{eqnarray*}
  \{(B(\hat\xi_+),\hat\xi_-)\in L^2(\sigma,\hat u_-^\ast TX)\times
  L^2(\sigma,\hat u_-^\ast TX)|\exists\,
  \xi_\pm\in \ker(D_{u_\pm}),\ \xi_\pm|_\sigma=\hat\xi_\pm\}
\end{eqnarray*}
is a Fredholm pair.

\begin{figure}[hpbt]
  \centering
  \includegraphics[width=7cm]{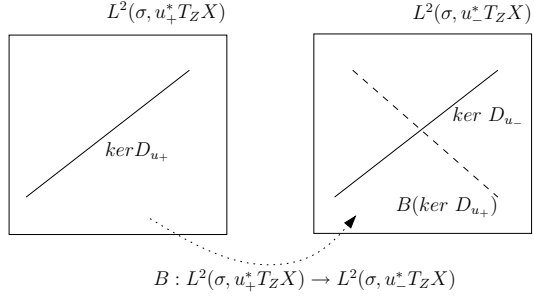}
  \caption{The map $B$ induces a Fredholm pair.}
  \label{fig:Fredholm_pair}
\end{figure}

To obtain nice compactifications it is important to first allow for
variation in the choice of parametrization of the closed
characteristic. This will induce an $S^1$-action on the space of
solutions, corresponding to the different choices of parametrization.
We refer to the quotient by the $S^1$ action as the {\em reduced
  moduli space}.
\begin{remark}{\em 
    We can generalize Theorem \ref{thm:fredholm} by incorporating
    variations of the folded domain and the choice of parametrization
    of the closed characteristic. Here are some brief comments on how
    this can be done.
    
    The first step is allowing variations in $j_0$ and $\tau_0$,
    modulo $Diff^+(\Sigma_0)$.  Note that $\tau_0$ is determined by
    the map by equation \ref{eq:pullback_det}. Since the folded
    diagonal is invariant under the gauge action of
    $\Diff^+(\Sigma_0)$, Theorem \ref{thm:fredholm} holds in this case
    with the index raised by the dimension of Teichm\"{u}ller space
    $-3\chi(\Sigma_0)$.
    
    Next consider variations in $\tau_1$, $j_1$, $\psi$ and $g$,
    modulo the action of the remaining factors of the gauge group
    $\Diff^+(\Sigma_1)\times \Map(\Sigma_1,\setR)$.  First note that
    for fixed domain location of the domain fold $\sigma_1$, the space
    of holomorphic diffeomorphisms
    $\psi^{-1}:\Sigma_0^+\rightarrow\Sigma_1^+$ sending $\sigma_0$ to
    $\sigma_1$ (with $j_0$ fixed and $j_1$ varying), has dimension
    $2\chi(\Sigma_1^+)+(1-3)\chi(\Sigma_1^+)=0$. Thus variations in
    $\psi$ and $j_1$ on $\Sigma_1^+$ do not change the dimension
    count.
    
    Finally, note that $g$ is determined by choice of $\tau_1$ and
    $\psi$ by equation (\ref{eq:folded domain condition}), and that
    the deformations of complex structure $j_1$ on $S$ are fixed by
    the tunneling map.  Thus the freedom left in choosing $\tau_1$,
    $j_1$, $\psi$,\ $g$ is exactly given by the remaining part
    $\Diff^+(\Sigma_1)\times\Map(\Sigma_1,\setR)$ of the gauge group
    $\mathcal{G}$. Changing the parametrization of the closed
    characteristic raises the index by 1. In conclusion, when varying
    the folded domain and the parametrization of the closed
    characteristic, and taking the quotient by the gauge group, the
    index is
    \begin{eqnarray}\label{eq:index2}
      \mu(u_+,K_+)+\mu(u_-,K_-)+(2-3)\chi(\Sigma)+1.
    \end{eqnarray}
  }
\end{remark}

\subsection{Homological Data}\label{sec:homological_data}
A folded map $(u_+,u_-)$ together with a pair of conjugate tunneling
maps $(v_+,v_-)$ gives rise to two relative homology classes
\begin{eqnarray*}
  A_\pm\in H_2(X^\pm,\mathcal R;\setZ),\qquad
  A_\pm=(u_\pm\sqcup_\sigma v_\pm)_\ast[\tilde\Sigma_\pm].
\end{eqnarray*}

Since the map gluing the tunneling domain $\hat S$ to $\Sigma_+$
($\Sigma_-$) is orientation preserving (reversing), and the
$\omega$-energies of the tunneling maps agree by Definition
(\ref{def:conjugate_tunneling}), we obtain the energy identities
\begin{eqnarray*}
  E_\omega(u_+)+E_\omega(v_+)&=&\omega\cdot A_+=\const\\
  E_\omega(u_-)-E_\omega(v_-)&=&\omega\cdot A_-=\const\\
  E_\omega(u_+)+E_\omega(u_-)&=&\omega\cdot(A_++A_-)=\const.
\end{eqnarray*}

Therefore the space of folded holomorphic maps breaks up into
components labeled by the relative homology classes $A_\pm$, and the
sum of the $\omega$-energies of the maps $u_+$ and $u_-$ is constant
in families.

\vskip1cm
\section{Examples of Folded Holomorphic Maps}
\label{sec:examples}
We give examples of folded holomorphic maps in two special cases.

\subsection{Folded Holomorphic Maps into Folded $E(1)$}
\label{sec:folded E(1)}
We come back to the example of $E(1)$ from Section
\ref{sec:motivating_example} and show that Definitions
\ref{def:folded_hol_E(1)} and \ref{def:folded_holomorphic_map}
coincide.

Note that in this case $d\alpha=0$ so equations (\ref{eq:conjugate_F})
and (\ref{eq:conjugate_L}) decouple. After fixing a closed
characteristic with parametrizations
$x_{\theta_0}(\theta)=(e^{2\pi i(\theta-\theta_0)},z_0)$ it is
straightforward to verify that the folded diagonal is given by the
graph of the family of functions
\begin{eqnarray*}
  \Phi_{\theta_0}:S^1\times T^2\rightarrow S^1\times T^2,\qquad
  \Phi_{\theta_0}(e^{2\pi i(\theta_0+\theta)},z)
  =(e^{2\pi i(\theta_0-\theta)},z).
\end{eqnarray*}
Note that the space of folded holomorphic maps of degree $d$ carries a
free $S^1$-action given by the choice of $\theta_0\in
S^1=\setR/\frac12\setZ$ since $\theta_0$ and
$\Phi_{\theta_0}=\Phi_{\theta_0+\frac12}$. We interpret this as the
different choices in gluing the folded symplectic $E(1)$.

This reproduces Definition \ref{def:folded_hol_E(1)}.

\subsection{Folded Holomorphic Rational Degree 1 Curves in $S^4$}
\label{sec:s4_example}
We explicitly characterize the moduli space of folded holomorphic
degree 1 rational curves by utilizing the symmetries of the folded
symplectic and complex structure on $S^4$ defined in Section
\ref{sec:motivating_example}. Essentially these curves come from
pseudo-holomorphic curves in $\setP^2$.

Recall that the pseudo-holomorphic cylinder over $S^3$
(``symplectization'') with its standard $\setR$-invariant structure is
biholomorphic to $\setC^2\setminus\{0\}$ via
\begin{eqnarray*}
  \Phi:\setR\times S^3\rightarrow\setC^2\setminus\{0\},
  \qquad(t,z)\mapsto e^{2t} z.
\end{eqnarray*}

For the rest of this section we fix homogeneous coordinates $[x:y:z]$
on $\setP^2$ and a corresponding embedding $\setC^2\subset \setP^2$,
$(z,w)\mapsto[x:w:1]$, whose complement is denoted by
$\setP^1_\infty$. Using this we can view finite asymptotic energy
pseudo-holomorphic maps in $\setR\times S^3$ as maps in $\setP^2$.
Conversely, pseudo-holomorphic maps in $\setP^2$ that have no
components that lie entirely in $\setP^1_\infty\cup\{0\}$ can be
viewed as (punctured) pseudo-holomorphic maps into $\setR\times S^3$
by restriction. A straightforward calculation reveals that punctured
finite asymptotic energy pseudo-holomorphic maps into $\setR\times
S^3$ extend over the punctures to pseudo-holomorphic maps into
$\setP^2$, and that conversely maps into $\setR\times S^3$ that are
restrictions of pseudo-holomorphic maps into $\setP^2$ have finite
asymptotic energy.

Note that
\begin{eqnarray*}
  H_2(B^4,\mathcal{R};\setZ)&=&S^2\times \setZ
\end{eqnarray*}
where the isomorphism is given by specifying the closed characteristic
and the multiplicity. We fix the closed characteristic parametrized by
the family
\begin{eqnarray*}
  x_m(\theta)=(m\,e^{2\pi\,i\theta},0)\subset S^3,\qquad 
  m\in S^1=\{z\in\setC|\,|z|=1\}
\end{eqnarray*}
and the {\em degree} $d=1$.

\begin{figure}[htbp]
  \centering
  \includegraphics[width=7cm]{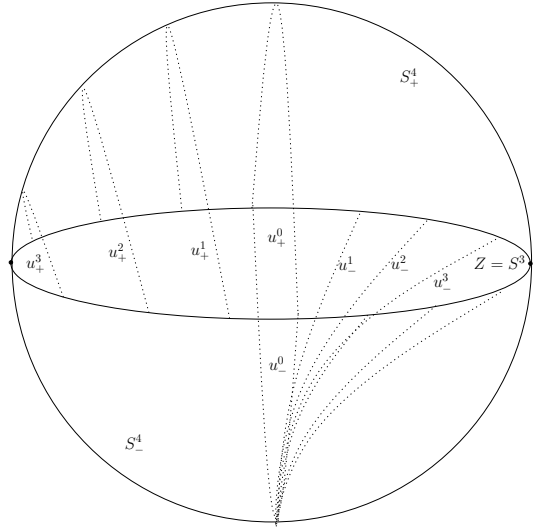}
  \caption{Several members of the moduli space of degree 1 maps into
    $S^4$. The map $(u_-^0,u_+^0)$ corresponds to the case $c=0$. The
    maps $u_+^n$ loose energy as they disappear into the fold whereas
    the maps $u_-^n$ gain energy.}
  \label{fig:S^4-Example}
\end{figure}

The open part of the moduli space is parametrized by
$c\in\setD=\{z\in\setC||z|<1\}$ and $m\in S^1$. Set
\begin{eqnarray*}
  \Sigma_+=\{z\in\hat\setC|\;|z|\le 1\}\qquad&
  \Sigma_-=\{z\in\hat\setC|\;|z|\ge 1\}&\\
  S=\{z\in\hat\setC|\;|z|\ge \sqrt{1-|c|^2}\}\qquad&
  \dot S=S\setminus\{\infty\}&
\end{eqnarray*}
where we work in $S^2=\hat\setC=\setC\cup\{\infty\}$. The isomorphism
$\psi:\Sigma_+\rightarrow\Sigma^+_1=\{z\in\setC|\,|z|\le\sqrt{1-|c|^2}\}$
is given by $z\mapsto \sqrt{1-|c|^2}z$. Then set
\begin{eqnarray*}
  u_+(z)=\sigma_+(\sqrt{1-|c|^2}m\,z,m\,c)\qquad&
  u_-(z)=\sigma_-(\sqrt{1-|c|^2}m/z,m\,c/z^2)&\\
  v_+(z)=\pi_{S^3}(m\,z,m\,c)\qquad&
  v_-(z)=\pi_{S^3}(m/z,m\,c/z^2).&
\end{eqnarray*}
We can parametrize this part of the moduli space, modulo
reparametrizations of the domain, by keeping track of one of the two
intersection points, say $(m\sqrt{1-|c|^2},mc)$, of $u_+$ and $u_-$,
which is in bijective correspondence with $Z\setminus y$, where $y$ is
the closed characteristic $y=\{(0,z)|\,|z|=1\}\in \mathcal{R}$.

\begin{figure}[hbtp]
  \centering
  \includegraphics[width=7cm]{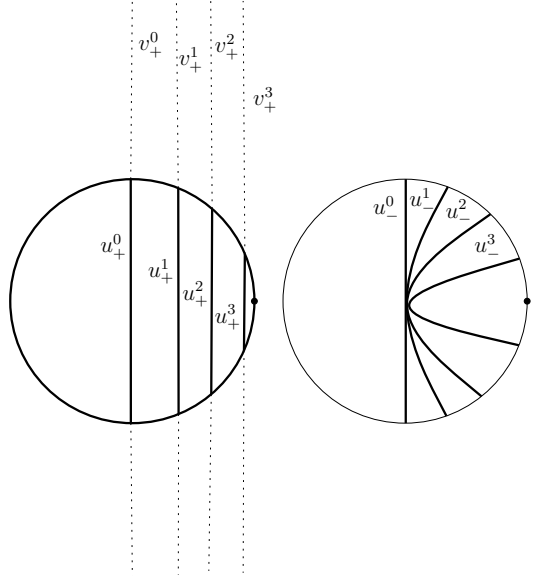}
  \caption{Here we visualize the degree 1 maps as maps into $\setC^2$. We
    also suspend the tunneling maps to ($\mathcal{H}$-) holomorphic
    maps into $\setC^2$.  The maps $u_\pm$ have image in the unit
    ball, whereas $v_+$ has image outside the unit ball. $v_-$
    coincides with $u_-$.
  }
  \label{fig:S^4-Example2}
\end{figure}

As the parameter $c$ leaves all compact subsets of $\setD$, the domain
$\Sigma_1$ degenerate so that $\dot S=\setC\setminus\{0\}$ is a sphere
with two punctures, so $v_\pm$ form a ``bubble'' in the fold. The
limit maps as $|c|\rightarrow 1$ $u_+(z)=(0,m\,c)\subset S^3$ sink
into the fold and become a point map with image on the closed
characteristic $y$. All energy is carried by $u_-(z)=(0,c/z^2)$, which
differ by a reparametrization of the domain for different choices of
$|c|=1$. So again we parametrize this portion of the moduli space,
modulo reparametrizations of the domain, by the corresponding
intersection point $(0,m\,c)\in y$, noting that any fixed choice of
$|c|=1$, while varying $m$, gives the same space of maps.

This shows that the moduli space of degree 1 rational maps has a
natural compactification to $S^3$. It carries a free $S^1$-action
given by the choice of parametrization $m$ of the closed
characteristic. The reduced moduli space, i.e. the quotient of the
moduli space by this action is naturally identified with $S^2$, with
quotient map the Hopf map.

\subsection{Rational Degree $d$ Maps into $S^4$}
\label{sec:degree_d}
We quickly indicate how to construct examples of rational degree $d$
maps into $S^4$. Again start from a degree $d$ curve $w:S^2\rightarrow
\setP^2$ that intersects the $\setP^1$ at infinity at the closed
characteristic $x$. Assume that $\Sigma_+=\{z\in S^2|\,||w(z)||\le
1\}$ in the homogeneous coordinates on the complement of the $\setP^1$
at infinity chosen above. Set $\Sigma_-=S^2\setminus \Sigma_+$ and
assume that $\Sigma_-$ is simply connected and $\pi_F\,dw\ne 0$ on
$\Sigma_-$. Then, with $S=\Sigma_-$ and the puncture corresponding to
the preimage of the $\setP^1$ at infinity we set
\begin{eqnarray*}
  u_+=\sigma_+\circ w|_{\Sigma_+}\qquad v_+=\pi_{S^3}w|_{\dot S}.
\end{eqnarray*}
To obtain the corresponding maps $v_-$ and $u_-$ let $\mu\in
S^1(T_pS)$ correspond to the direction given by the positive real line
in $\setC$ and let $f:\dot S\rightarrow \setC$ be the unique
holomorphic function satisfying $x_m(\Im(f(p)))=-2v_+(\mu)$, where we
use the group structure induced by $x_m$, and
$d[\Re(f)]=-2v_+^\ast\alpha\circ j$ on $T\sigma$ and set
\begin{eqnarray*}
  u_-(z)=\sigma_-(f(z)\cdot w|_{\Sigma_-}(z))\qquad
  v_-=\pi_{S^3}u_-|_{\dot S}.
\end{eqnarray*}

To compute the dimension of the moduli space of degree $d$ maps we
compute the index of the linearized operator (\ref{eq:index2}) at a
map satisfying $\pi_F\,du_\pm\ne 0$ along the domain fold $\sigma$. By
assumption the domains $\Sigma_\pm$ are both diffeomorphic to the disk
$\setD$. Homotope the boundary conditions $F_\pm(z)$ to $\setR\cdot
z^d\times \setR\subset u_\pm^\ast(\setC\times \setC)$, so
\begin{eqnarray*}
  \mu(H_\pm,F_\pm)=\mu(\setC\times\setC,\setR\cdot z^d\times\setR)
  =\mu(\setC,\setR\cdot z^d)+\mu(\setC,\setR)=2d.
\end{eqnarray*}
Therefore
\begin{eqnarray*}
  index=\mu(H_+,F_+)+\mu(H_-,F_-)+(2-3)\chi(S^2)+1=4d-1.
\end{eqnarray*}
Again, these spaces carry a free $S^1$ action corresponding to the
choice $m\in S^1$ of parametrization of the closed
characteristic.

\vskip1cm
\section{Tunneling Maps in Symplectic Manifolds}
\label{sec:symplectic_tunneling}
Here we want to briefly explain how tunneling maps come up in the
usual symplectic setting. From this point of view, tunneling maps
appear as tools for studying $J$-holomorphic curves relative to a
codimension 1 hypersurface in a symplectic manifold. We transfer our
definitions to this case:

Let $(X,\omega)$ be a symplectic manifold and $f:X\rightarrow \setR$ a
smooth function with transverse zeros. Then $Z=f^{-1}$ is a smooth
hypersurface, separating $X$ into two parts labeled $X_+$ and $X_-$ by
the sign of $f$ on them. Assume that there exists a 1-form $\alpha$ on
$Z$ so that $Z$ together with $\alpha$ and $\omega$ admits an
$S^1$-invariant structure as in Definition
\ref{def:S^1-invariant_fold}. Choose a compatible almost complex
structure $J$ on $X$ that is $S^1$-invariant over $Z$.

Fix a folded domain as in Definition \ref{def:folded domain} and let
$u_\pm:\Sigma_\pm\rightarrow X_\pm$ be $J$-holomorphic with $\tau=u^\ast
f$. Assume $\tau$ vanishes transversely and
$\sigma=\tau^{-1}(0)\neq\emptyset$, so $\sigma$ is a smooth non-empty
compact submanifold of $\Sigma$, separating $\Sigma$.

Now, just like in the folded symplectic case we may look for tunneling
maps in $Z$ that connect the image of $u|_\sigma$ to closed
characteristics, i.e. an $\mathcal{H}$-holomorphic map $v:\dot
S\rightarrow Z$ with $v|_\sigma=u|_\sigma$. One might at first expect
that tunneling maps should be pseudo-holomorphic instead of
$\mathcal{H}$-holomorphic, but note that the index for
pseudo-holomorphic tunneling maps that are immersions is negative the
first Betty number of the domain $-b^1(S)$.

As opposed to the folded symplectic case, the complex structures
$J_\pm$ induced on $T_ZX$ coming from $X_+$ and $X_-$ agree in the
symplectic setting. Thus the argument in the discussion of the sign of
$u^\ast\alpha$ in Remark \ref{rem:orientations} has to be modified.
Therefore the folded diagonal $\Delta^Z$ does not pose Fredholm
boundary conditions in the symplectic case as becomes clear in the
proof of Theorem \ref{thm:fredholm}.

To understand this better we take another look at tunneling maps in
the folded symplectic setting. Let $(v_+,v_-)$ be conjugate tunneling
maps and consider the suspension $\tilde v_\pm$ of tunneling maps
$v_\pm$ into $\setR\times Z$, where $v_\pm$ is a
$\mathcal{H}$-holomorphic map with respect to $J_\pm$.  Assume for
simplicity that $S$ is a punctured disk, so the tunneling maps are
actually $J_\pm$-holomorphic. Then the equation (\ref{eq:conjugate_L})
of the definition of conjugate tunneling maps can be rewritten in
terms of the $\setR$-components $a_+$ and $a_-$ of the corresponding
tunneling maps. Then
\begin{eqnarray*}
  \tilde v_+^\ast\alpha\circ j+\tilde v_-^\ast\alpha\circ j
  =da_+-da_-
\end{eqnarray*}
and using the fact that the $\setR$-component is only determined up to
a constant, equation (\ref{eq:conjugate_L}) says
\begin{eqnarray*}
  a_+=a_-\qquad \mathrm{on}\ \sigma.
\end{eqnarray*}
Written in this way the equation carries over verbatim to the
(non-folded) symplectic setting.

We can follow the above transformations backward and obtain the
replacement of equation (\ref{eq:conjugate_L}) for the symplectic
setting:
\begin{eqnarray*}
  v_+^\ast\alpha\circ j=v_-^\ast\alpha\circ j\ \mathrm{on}\ T\sigma
\end{eqnarray*}
for $\mathcal{H}$-holomorphic maps $v_\pm$.

But this, together with the remaining equations in Definition
\ref{def:conjugate_tunneling} implies that $v_+=v_-$. Thus in the
symplectic case, the analogue of the folded diagonal is the actual
diagonal 
\begin{eqnarray*} 
  \Delta=\{(\hat v,\hat v)|\hat v:\sigma\rightarrow Z\}\subset
  \Map(\sigma,Z)\times \Map(\sigma,Z).
\end{eqnarray*}
Viewing this the other way, the folded diagonal $\Delta^Z$ in the
folded symplectic setting is analogue of the actual diagonal $\Delta$
in the symplectic setting.

\vspace{2cm}
\noindent{\sc\bf\large Acknowledgments:}\\
I would like to thank my adviser Thomas H. Parker for his support,
guidance and helpful discussions.

\bibliography{folded}
\end{document}